# Kullback Leibler property of kernel mixture priors in Bayesian density estimation

## Yuefeng Wu

*Department of Statistics, North Carolina State University*
*e-mail:* ywu9@ncsu.edu

## Subhashis Ghosal

*Department of Statistics, North Carolina State University*
*e-mail:* sghosal@stat.ncsu.edu

**Abstract:** Positivity of the prior probability of Kullback-Leibler neighborhood around the true density, commonly known as the Kullback-Leibler property, plays a fundamental role in posterior consistency. A popular prior for Bayesian estimation is given by a Dirichlet mixture, where the kernels are chosen depending on the sample space and the class of densities to be estimated. The Kullback-Leibler property of the Dirichlet mixture prior has been shown for some special kernels like the normal density or Bernstein polynomial, under appropriate conditions. In this paper, we obtain easily verifiable sufficient conditions, under which a prior obtained by mixing a general kernel possesses the Kullback-Leibler property. We study a wide variety of kernel used in practice, including the normal, $t$, histogram, gamma, Weibull densities and so on, and show that the Kullback-Leibler property holds if some easily verifiable conditions are satisfied at the true density. This gives a catalog of conditions required for the Kullback-Leibler property, which can be readily used in applications.

## 1. Introduction

Density estimation, which is also relevant in various applications such as cluster analysis and robust estimation, is a fundamental nonparametric inference problem. In Bayesian approach to density estimation, a prior such as a Gaussian process, a Polya tree process, or a Dirichlet mixture is constructed on the space of probability densities. Dirichlet mixtures were introduced by Ferguson [9] and Lo [21] who also obtained expressions for resulting posterior and predictive distribution. West [30], West, Müller and Escobar [31] and Escobar and West [6; 7] developed powerful Markov chain Monte Carlo methods to calculate Bayes estimates and other posterior quantities for Dirichlet mixtures.

The priors of interest in this paper are of mixture type and can be described in terms of a kernel and a prior for the mixing distribution. Let $\mathfrak{X}$ be the sample space and $\Theta$ is the space of the mixing parameter $\theta$. Let $K(x; \theta)$ be the kernel on $\mathfrak{X} \times \Theta$, i.e., $K(x; \theta)$ is a jointly measurable function such that for all $\theta$, $K(\cdot; \theta)$ is a probability density on $\mathfrak{X}$. The choice of an appropriate kernel depends on the underlying sample space $\mathfrak{X}$, on which the true density is defined. If $\mathfrak{X}$ is the entire real line, a location-scale kernel is appropriate. If $\mathfrak{X}$ is the unit interval, a uniform or triangular density kernel, or Bernstein polynomial may be considered. If $\mathfrak{X}$ is the positive half line $(0, \infty)$, mixtures of gamma, Weibull, lognormal, exponential or inverse gamma may be used. Petrone and Veronese [25] discussed the issue of the choice of a kernel in view of a constructive approximation known as the Feller sampling scheme. Let $P$, the mixing distribution on $\Theta$, be given a prior $\Pi$ on $\mathscr{M}(\Theta)$, the space of probability measure on $\Theta$. Let supp($\Pi$) denote the weak support of $\Pi$. The prior on $P$ and the chosen kernel then give rise to a prior on $\mathscr{D}(\mathfrak{X})$, the space of densities on $\mathfrak{X}$, via the map $P \mapsto f_P(x) := \int K(x; \theta) dP(\theta)$. We shall call such a prior a type I mixture prior or Prior 1 in short. To enrich the family of the kernels, let the kernel function contain another parameter $\phi$, referred to as the hyper parameter. In this case, we shall denote the kernel by $K(x; \theta, \phi)$. The hyper parameter $\phi$ might be elicited a priori or be given a prior. In the former case, such a prior essentially reduces to Prior 1. For the latter case, assume that $\phi$ is independent of $P$ and denote the prior for $\phi$ by $\mu$. Let $\Phi$ be the space of $\phi$ and supp($\mu$) denote the support of $\mu$. With such a random hyper parameter in the chosen kernel, the prior on densities is induced by $\mu \times \Pi$ via the map $(\phi, P) \mapsto f_{P,\phi}(x) := \int K(x; \theta, \phi) dP(\theta)$. We shall call this prior a Type II mixture prior or simply Prior 2. Clearly, Prior 2 contains Prior 1 as a special case where $\phi$ is treated as a vacuous parameter. In some situations, the prior $\Pi$ may contain an additional indexing parameter $\xi$. For instance, when $\Pi$ is the Dirichlet process with base measure $\alpha_\xi$ (written as DP($\alpha_\xi$)) depending on an indexing parameter $\xi$, which is also given a prior, we obtain a mixture of Dirichlet processes (MDP) [1] prior for mixing distribution $P$. Addition of this hierarchical structure to Prior 1 or Prior 2 gives somewhat more flexibility. In this paper, we do not make any specific assumption on $\Pi$ like DP or MDP other than requiring that it has large weak support. The prior induced on the space of densities by a mixing distribution $P \sim \Pi$ (and $\phi \sim \mu$ and $\xi \sim \pi$) will be denoted by $\Pi^*$ and we shall refer to it as a kernel mixture prior. Note that the variable $x$ and the parameters $\theta$, $\phi$ and $\xi$ mentioned above are not necessarily one-dimensional.

Asymptotic properties, such as consistency, and rate of convergence of the posterior distribution based on kernel mixture priors were established by Ghosal, Ghosh and Ramamoorthi [11], Tokdar [29], and Ghosal and van der Vaart [13; 14], when the kernel is chosen to be a normal probability density (and the prior distribution of the mixing distribution is DP). Similar results for Dirichlet mixture of Bernstein polynomials were shown by Petrone and Wasserman [26], Ghosal [10] and Kruijer and van der Vaart [19]. However, in the literature, there is a lack of such results for mixture of other kernels, which are also widely used in practice. We are only aware of the article by Petrone and Veronese [25] who considered general kernels. However, they derived consistency only under the strong and unrealistic condition that the true density is exactly of the mixture type for some compactly supported mixing distribution, or the true density itself is compactly supported and is approximated in terms of Kullback-Leibler divergence by its convolution with the chosen kernel.

Schwartz [28] showed that the consistency at a true density $f_0$ holds if the prior assigns positive probabilities to specific type of neighborhoods of $f_0$ defined by Kullback-Leibler divergence measure and the size of the model is restricted in some appropriate sense. Thus the prior positivity condition, known as the Kullback-Leibler property (KL property), is fundamental in posterior consistency studies. More formally, let a density function $f$ be given a prior $\Pi^*$. Define a Kullbck-Leibler neighborhood of $f$ of size $\epsilon$ by $\mathscr{K}_\epsilon(f) = \{g : \mathcal{K}(f; g) < \epsilon\}$, where $\mathcal{K}(f; g) = \int f \log(f/g)$, the Kullback-Leibler divergence between $f$ and $g$. We say that the KL property holds at $f_0 \in \mathscr{D}(\mathfrak{X})$ or $f_0$ is in the Kullback-Leibler support (KL support) of $\Pi^*$, and write $f_0 \in \mathrm{KL}(\Pi^*)$, if $\Pi^*(\mathscr{K}_\epsilon(f_0)) > 0$ for every $\epsilon > 0$. For the weak topology, the size condition in Schwartz's theorem holds automatically [16, Theorem 4.4.2]. Further, Ghosal, Ghosh and Rammamoorthi [12] argued that this property drives consistency of the parametric part in some semiparametric models.

This paper addresses issues about KL property of general kernel mixture priors, thus addressing one of the most important issues in posterior consistency. We discuss the KL property for general kernel mixture priors, which are not restricted by any particular type of kernel or by a prior distribution for mixing distribution. The distinguished feature of our results is that we allow the true density to be not of the chosen mixture type, and impose only simple moment conditions and qualitative conditions like continuity or positivity.

Ghosal, Ghosh and Rammamoorthi [11] presented results on consistency for Dirichlet location mixture of a normal kernel with an additional scale parameter in terms of both weak and $L_1$-topologies. Tokdar [29, Theorem 3.2] considered a location-scale mixture of the normal kernel and established consistency in weak topology (weak consistency) under more relaxed conditions. If the prior $\Pi$ is chosen to be $\mathrm{DP}(\alpha)$, Tokdar [29] also weakened a moment condition on the true density in his Theorem 3.3. His Theorem 3.2 will be implied by Theorem 4 in this paper (with the choice $\lambda = 0$ there). In fact, we establish the KL property for a general location-scale kernel mixture and show that such a result applies to various kernels including the skew-normal, $t$, double-exponential and logistic. This is a substantial generalization of results known for only the normal kernel thus far. Moreover, we obtain results about the KL property for priors with kernels not belonging to location–scale families, e.g., the Weibull, gamma, uniform, and exponential kernels. The examples studied here provide a ready catalog of conditions required for the KL property to hold for virtually all kernel mixture priors that are of practical interest.

With the the help of our results on KL property, consistency in $L_1$- (equivalently, Hellinger) distance can be obtained by constructing appropriate sieves approximating the class of mixtures and establishing entropy bounds for them. Since the techniques used for sieve construction and bounding entropy vary widely depending on the chosen kernel, we do not address $L_1$-consistency in this paper.

The paper is organized as follows. In Section 2, we study the kernel mixture priors under complete generality without specifying a kernel or the nature of it. In Section 3, using the results provided in Section 2, we study the priors with kernels of the location-scale type. In Section 4, the priors with concretely specified kernels are studied as examples by using the results in the previous sections.

## 2. General Kernel Mixture Priors

First we observe that the Kullback Leibler property is preserved under taking mixtures.

**Lemma 1.** *Let $f|\xi \sim \Pi_\xi^*$, where $\xi$ is an indexing parameter following a prior $\pi$ and let $f_0$ be the true density. Suppose that there exists a set $B$ with properties $\Pi(B) > 0$ and $B \subset \{\xi : f_0 \in \mathrm{KL}(\Pi_\xi^*)\}$. Then $f_0 \in \mathrm{KL}(\Pi^*)$, where $\Pi^* = \int \Pi_\xi^* d\pi(\xi)$.*

The proof is almost a trivial application of Fubini's theorem, since

$$\Pi^*(f : \mathcal{K}(f_0; f) < \epsilon) \geq \int_B \Pi_\xi^*(f : \mathcal{K}(f_0; f) < \epsilon) d\pi(\xi) > 0.$$

In view of this result, henceforth we shall discard the indexing parameter $\xi$ from our prior.

**Theorem 1.** *Let $f_0$ be the true density, $\mu$ and $\Pi$ be priors for the hyper parameter and the mixing distribution in Prior 2, and $\Pi^*$ be the prior induced by $\mu$ and $\Pi$ on $\mathscr{D}(\mathfrak{X})$. If for any $\epsilon > 0$, there exists $P_\epsilon$, $\phi_\epsilon$, $A \subset \Phi$ with $\mu(A) > 0$ and $\mathscr{W} \subset \mathscr{M}(\Theta)$ with $\Pi(\mathscr{W}) > 0$, such that*

*A1. $\int f_0 \log \frac{f_0}{f_{P_\epsilon, \phi_\epsilon}} < \epsilon$,*

*A2. $\int f_0 \log \frac{f_{P_\epsilon, \phi_\epsilon}}{f_{P_\epsilon, \phi}} < \epsilon$ for every $\phi \in A$, and*

*A3. $\int f_0 \log \frac{f_{P_\epsilon, \phi}}{f_{P, \phi}} < \epsilon$ for every $P \in \mathscr{W}$, $\phi \in A$,*

*then $f_0 \in \mathrm{KL}(\Pi^*)$.*

PROOF. For any $\epsilon > 0$, $\phi \in A$ and $P \in \mathscr{W}$,

$$
\begin{aligned}
\int_{\mathfrak{X}} f_0(x) \log \frac{f_0(x)}{f_{P,\phi}(x)} dx &= \int_{\mathfrak{X}} f_0(x) \log \frac{f_0(x)}{f_{P_\epsilon, \phi_\epsilon}(x)} dx \\
&\quad + \int_{\mathfrak{X}} f_0(x) \log \frac{f_{P_\epsilon, \phi_\epsilon}(x)}{f_{P_\epsilon, \phi}(x)} dx \\
&\quad + \int_{\mathfrak{X}} f_0(x) \log \frac{f_{P_\epsilon, \phi}(x)}{f_{P,\phi}(x)} dx < 3\epsilon, \quad (1)
\end{aligned}
$$

Hence,

$$\Pi^*\{f : f \in \mathscr{K}_{3\epsilon}(f_0)\} \geq \Pi^*\{f_{P,\phi} : P \in \mathscr{W}, \phi \in A\} = (\Pi \times \mu)(\mathscr{W} \times A) > 0.$$

$\square$

**Remark 1.** If $\Pi = \mathrm{DP}(\alpha)$ and $\mathrm{supp}(P_\epsilon) \subset \mathrm{supp}(\alpha)$, then $P_\epsilon \in \mathrm{supp}(\Pi)$; see, for instance, Theorem 3.2.4 of [16]. In particular, the condition holds for any chosen $P_\epsilon$ if $\alpha$ is fully supported on $\Theta$. A similar assertion holds when $\Pi$ is the Polya tree prior $\mathrm{PT}(\{\mathcal{T}_m\}, \mathscr{A})$ (see [20]). Let $\mathcal{T}_m$ be a collection of gradually refining binary partitions and $\mathscr{A} = \{\alpha_{\epsilon_1,\ldots,\epsilon_m} : \epsilon_1,\ldots,\epsilon_m = 0 \text{ or } 1, m \geq 1\}$. If the end points of $\mathcal{T}_m$ form a dense subset of some set $\mathscr{S}$ where $\mathscr{S} \supset \mathrm{supp}(P_\epsilon)$ and the elements of $\mathscr{A}$, which control the beta distributions regulating the mass allocation to the sets in $\Pi_m$, are positive, then also $P_\epsilon \in \mathrm{supp}(\Pi)$. This is implicit in Theorem 5 of [20] or Theorem 3.3.6 of [20]; for an explicit statement and proof, see Theorem 2.20 of [15]. Now, if $\mathscr{W}$ is an open neighborhood of $P_\epsilon$, then $\Pi(\mathscr{W}) > 0$ holds.

**Remark 2.** Assume that $\phi_\epsilon \in \mathrm{supp}(\mu)$. Condition A2 clearly holds with $A$ an open neighborhood of $\phi_\epsilon$, assuming that $\phi \mapsto \int f_0 \log(f_{P_\epsilon,\phi_\epsilon}/f_{P_\epsilon,\phi})$ is continuous.

In most application, we can choose $P_\epsilon$ to be compactly supported. Compactness of $\mathrm{supp}(P_\epsilon)$ often helps satisfy condition A4–A9 in Lemmas 2 and 3, which are useful in verifying the conditions of Theorem 1.

**Lemma 2.** *Let $f_0$, $\Pi$, $\mu$ and $\Pi^*$ be the same as in Theorem 1. If for any $\epsilon > 0$, there exist $P_\epsilon$, a set $D$ containing $\mathrm{supp}(P_\epsilon)$, and $\phi_\epsilon \in \mathrm{supp}(\mu)$ such that A1 holds and the kernel function $K$ satisfies*

*A4. for any given $x$ and $\theta$, the map $\phi \mapsto K(x;\theta,\phi)$ is continuous on the interior of the support of $\mu$;*

*A5. $\int_{\mathfrak{X}} \left\{ \left| \log \frac{\sup_{\theta \in D} K(x;\theta,\phi_\epsilon)}{\inf_{\theta \in D} K(x;\theta,\phi)} \right| + \left| \log \frac{\sup_{\theta \in D} K(x;\theta,\phi)}{\inf_{\theta \in D} K(x;\theta,\phi_\epsilon)} \right| \right\} f_0(x)dx < \infty$ for every $\phi \in N(\phi_\epsilon)$, where $N(\phi_\epsilon)$ is an open neighborhood of $\phi_\epsilon$;*

*A6. for any given $x \in \mathfrak{X}$, $\theta \in D$ and $\phi \in N(\phi_\epsilon)$, there exists $g(x,\theta)$ such that $g(x,\theta) \geq K(x;\theta,\phi)$, and $\int g(x,\theta)dP_\epsilon(\theta) < \infty$;*

*then there exists a set $A \subset \Phi$ such that A2 holds.*

Proof. By Condition A4, we have that $K(x;\theta,\phi) \to K(x;\theta,\phi_\epsilon)$ as $\phi \to \phi_\epsilon$, for any given $x$ and $\theta$. By Condition A6 and the dominated convergence theorem (DCT), $f_{P_\epsilon,\phi}(x) \to f_{P_\epsilon,\phi_\epsilon}(x,)$ as $\phi \to \phi_\epsilon$, for any given $x$. Equivalently, this can be written as

$$\log \frac{f_{P_\epsilon,\phi_\epsilon}}{f_{P_\epsilon,\phi}} \to 0 \text{ pointwise, as } \phi \to \phi_\epsilon. \tag{2}$$

Note that

$$\left| \log \frac{f_{P_\epsilon,\phi_\epsilon}}{f_{P_\epsilon,\phi}} \right| \leq \begin{cases} \left| \log \frac{\sup_{\theta \in D} K(x;\theta,\phi_\epsilon)}{\inf_{\theta \in D} K(x;\theta,\phi)} \right|, & \text{if } \frac{f_{P,\phi_\epsilon}}{f_{P,\phi}} \geq 1, \\[2mm] \left| \log \frac{\sup_{\theta \in D} K(x;\theta,\phi)}{\inf_{\theta \in D} K(x;\theta,\phi_\epsilon)} \right|, & \text{if } \frac{f_{P,\phi_\epsilon}}{f_{P,\phi}} < 1. \end{cases}$$

By Condition A5 and the DCT, $\int f_0 \log \frac{f_{P_\epsilon,\phi_\epsilon}}{f_{P_\epsilon,\phi}} \to 0$ as $\phi \to \phi_\epsilon$. Hence, for given $\epsilon > 0$, there exists $\delta > 0$ such that $\int f_0 \log \frac{f_{P_\epsilon,\phi_\epsilon}}{f_{P_\epsilon,\phi}} < \epsilon$ if $|\phi - \phi_\epsilon| < \delta$. If $A = \{\phi : |\phi - \phi_\epsilon| < \delta\} \cap N(\phi_\epsilon)$, then $\int f_0 \log \frac{f_{P_\epsilon, \phi_\epsilon}}{f_{P_\epsilon, \phi}} < \epsilon$ for all $\phi \in A$. The proof is completed by noticing that $\mu(A) > 0$, since $A$ is an open neighborhood of $\phi_\epsilon \in \text{supp}(\mu)$. □

**Lemma 3.** *Let $f_0$, $\Pi$, $\mu$ and $\Pi^*$ be the same as in Theorem 1. If for any $\epsilon > 0$, there exist $P_\epsilon \in \text{supp}(\Pi)$, $\phi_\epsilon \in \text{supp}(\mu)$, and $A \subset \Phi$ with $\mu(A) > 0$ such that Conditions A1 and A2 hold and for some closed $D \supset \text{supp}(P_\epsilon)$, the kernel function $K$ and prior $\Pi$ satisfy*

*A7. for any $\phi \in A$, $\int \log \frac{f_{P_\epsilon, \phi}(x)}{\inf_{\theta \in D} K(x, \theta, \phi)} f_0(x) dx < \infty$;*
*A8. $c := \inf_{x \in C} \inf_{\theta \in D} K(x; \theta, \phi) > 0$, for any compact $C \subset \mathfrak{X}$;*
*A9. for any given $\phi \in A$ and compact $C \subset \mathfrak{X}$, there exists $E$ containing $D$ in its interior such that the family of maps $\{\theta \mapsto K(x; \theta, \phi), x \in C\}$ is uniformly equicontinuous on $E \subset \Theta$, and $\sup\{K(x; \theta, \phi) : x \in C, \theta \in E^c\} < c\epsilon/4$;*

*then there exists $\mathcal{W} \subset \mathcal{M}(\Theta)$ such that Condition A3 holds and $\Pi(\mathcal{W}) > 0$.*

PROOF. For any $\phi \in A$, write

$$
\begin{aligned}
\int_{\mathfrak{X}} f_0(x) \log \frac{f_{P_\epsilon, \phi}(x)}{f_{P, \phi}(x)} dx &= \int_{C^c} f_0(x) \log \frac{f_{P_\epsilon, \phi}(x)}{f_{P, \phi}(x)} dx \\
&\quad + \int_C f_0(x) \log \frac{f_{P_\epsilon, \phi}(x)}{f_{P, \phi}(x)} dx.
\end{aligned}
\tag{3}
$$

Now, since $P_\epsilon(D) = 1 > \frac{1}{2}$, $\mathcal{V} = \{P : P(D) > \frac{1}{2}\}$ is an open neighborhood of $P_\epsilon$ by the Portmanteau Theorem. For any $P \in \mathcal{V}$ and $\phi \in A$,

$$
\begin{aligned}
&\int_{C^c} f_0(x) \log \frac{f_{P_\epsilon, \phi}(x)}{f_{P, \phi}(x)} dx \\
&\leq \int_{C^c} f_0(x) \log \frac{f_{P_\epsilon, \phi}(x)}{\int_{\theta \in D} \inf_{\theta \in D} K(x; \theta, \phi) dP(\theta)} dx \\
&\leq \int_{C^c} f_0(x) \log \frac{f_{P_\epsilon, \phi}(x)}{\inf_{\theta \in D} K(x; \theta, \phi)} dx + (\log 2) P_{f_0}(C^c);
\end{aligned}
$$

here $P_{f_0}$ is the probability measure corresponding to $f_0$. By Condition A7, there exists compact $C \subset \mathfrak{X}$, such that

$$
\int_{C^c} f_0(x) \log \frac{f_{P_\epsilon, \phi}(x)}{\inf_{\theta \in D} K(x; \theta, \phi)} dx < \epsilon/4.
\tag{4}
$$

We can further ensure that $P_{f_0}(C^c) < \epsilon/4$, so the bound for $\int_{C^c} f_0 \log \frac{f_{P_\epsilon, \phi}}{f_{P, \phi}}$ is less than $\epsilon/2$. Now, if we can show that for the given $\epsilon > 0$, there exists a weak neighborhood $\mathcal{U}$ of $P_\epsilon$, such that $\int_C f_0(x) \log \frac{f_{P_\epsilon, \phi}(x)}{f_{P, \phi}(x)} dx < \epsilon/2$ for any $P \in \mathcal{U}$ and $\phi \in A$, then Lemma 3 is proved by letting $\mathcal{W} = \mathcal{U} \cap \mathcal{V}$.

Observing that for any given $\phi \in A$, the family of maps $\{\theta \mapsto K(x; \theta, \phi) : x \in C\}$ is uniformly equicontinuous on $E \subset \Theta$, by the Arzela-Ascoli theorem,

(see [27, pp. 169]) for any $\delta > 0$, there exist $x_1, x_2, \ldots, x_m$, such that, for any $x \in C$,

$$\sup_{\theta \in E} |K(x; \theta, \phi) - K(x_i; \theta, \phi)| < c\delta. \tag{5}$$

for some $i = 1, 2, \ldots, m$.

Let $\mathscr{U} = \{P : |\int_E K(x_i; \theta, \phi) dP_\epsilon(\theta) - \int_E K(x_i; \theta, \phi) dP(\theta)| < c\delta, \quad i = 1, 2, \ldots, m\}$. Then $\mathscr{U}$ is an open weak neighborhoods of $P_\epsilon$ since $P_\epsilon \in \text{supp}(\Pi)$ and $P_\epsilon(\partial E) = 0$. For any $x \in C$, choosing $x_i$ to satisfy (5), we have that

$$\begin{aligned}
&\left| \int_\Theta K(x; \theta, \phi) dP(\theta) - \int_\Theta K(x; \theta, \phi) dP_\epsilon(\theta) \right| \\
&\quad \leq \sup\{K(x; \theta, \phi) : \theta \in E^c, x \in C\} \\
&\qquad + \left| \int_E K(x; \theta, \phi) dP(\theta) - \int_E K(x_i; \theta, \phi) dP(\theta) \right| \\
&\qquad + \left| \int_E K(x_i; \theta, \phi) dP(\theta) - \int_E K(x_i; \theta, \phi) dP_\epsilon(\theta) \right| \\
&\qquad + \left| \int_E K(x_i; \theta, \phi) dP_\epsilon(\theta) - \int_E K(x; \theta, \phi) dP_\epsilon(\theta) \right| \\
&\quad < \frac{c\epsilon}{4} + 2c\delta + \left| \int_E K(x_i; \theta, \phi) dP_\epsilon(\theta) - \int_E K(x_i; \theta, \phi) dP(\theta) \right| \\
&\quad < c\left( \frac{\epsilon}{4} + 3\delta \right)
\end{aligned} \tag{6}$$

if $P \in \mathscr{U}$. Also $\int_\Theta K(x; \theta, \phi) dP_\epsilon(\theta) > c$ for any $x \in C$, since $P_\epsilon$ has support in $D$. Hence, given $\phi \in A$, for any $P \in \mathscr{U}$ and $x \in C$,

$$\left| \frac{\int_\Theta K(x; \theta, \phi) dP(\theta)}{\int_\Theta K(x; \theta, \phi) dP_\epsilon(\theta)} - 1 \right| \leq 3\delta + \frac{\epsilon}{4}.$$

Then, for $3\delta + \epsilon/4 < 1$,

$$\left| \frac{\int_\Theta K(x; \theta, \phi) dP_\epsilon(\theta)}{\int_\Theta K(x; \theta, \phi) dP(\theta)} - 1 \right| < \frac{3\delta + \epsilon/4}{1 - 3\delta - \epsilon/4}.$$

By choosing $\delta$ small enough, we can ensure that the right hand side (RHS) of the last display is less than $\epsilon/2$. Hence, for any given $\phi \in A$

$$\int_C f_0(x) \log \frac{f_{P_\epsilon, \phi}(x)}{f_{P, \phi}(x)} dx \leq \sup_{x \in C} \left| \frac{\int_\Theta K(x; \theta, \phi) dP_\epsilon(\theta)}{\int_\Theta K(x; \theta, \phi) dP(\theta)} - 1 \right| < \epsilon/2$$

for any $P \in \mathscr{U}$. $\qquad \square$

## 3. Location scale kernel

In this section we discuss priors with kernel functions belonging to location scale families. We write the kernels as $K(x; \theta, h) = \frac{1}{h^d} \chi(\frac{x-\theta}{h})$, where $\chi(\cdot)$ is a probability density function defined on $\mathbb{R}^d$, $x = (x_1, \ldots, x_d)$, and $\theta = (\theta_1, \ldots, \theta_d)$ are $d$-dimensional vectors and $h \in (0, \infty)$. Let $\|x\|$ denote $\sqrt{x_1^2 + x_2^2 + \ldots + x_d^2}$, and $\chi_i'(x)$ denote $\frac{\partial \chi(x)}{\partial x_i}$. Obviously, when $d = 1$, this reduces to ordinary derivative and $\|\cdot\|$ denotes absolute value. We have the following theorems, whose proofs use some ideas from the proof of Theorem 3.2 of [29].

**Theorem 2.** *Let $f_0(x)$ be the true density and $\Pi^*$ be a type I prior on $\mathscr{D}(\mathfrak{X})$ with kernel function $h^{-d}\chi(\frac{x-\theta}{h})$, i.e. $P \sim \Pi$, and given $P$, $(\theta, h) \sim P$. If $\chi(\cdot)$ and $f_0(x)$ satisfy:*

*B1. $\chi(\cdot)$ is bounded, continuous and positive everywhere;*

*B2. there exists $l_1 > 0$ such that $\chi(x)$ decreases as $x$ moves away from $0$ outside the ball $\{x : \|x\| < l_1\}$;*

*B3. there exists $l_2 > 0$ such that $\sum_{i=1}^{d} z_i \frac{\chi_i'(z)}{\chi(z)} < -1$ for $\|z\| \geq l_2$ and $i = 1, \ldots, d$;*

*B4. for some $0 < M < \infty$, $0 < f_0(x) \leq M$ for all $x$;*

*B5. $|\int f_0(x) \log f_0(x) dx| < \infty$;*

*B6. for some $\delta > 0$, $\int f_0(x) \log \frac{f_0(x)}{\phi_\delta(x)} dx < \infty$, where $\phi_\delta(x) = \inf_{\|t-x\| < \delta} f_0(t)$;*

*B7. there exists $\eta > 0$, such that $|\int f_0(x) \log \chi(2x\|x\|^\eta) dx| < \infty$ and $\int f_0(x)|\log \chi(\frac{x-a}{b})| dx < \infty$ for any $a \in \mathbb{R}^d$, $b \in (0, \infty)$;*

*B8. the weak support of $\Pi$ is $\mathscr{M}(\mathbb{R}^d \times \mathbb{R}^+)$;*

*B9. when $d \geq 2$, $\chi(y) = o(\|y\|^{-d})$ as $\|y\| \to \infty$.*

*Then $f_0 \in \mathrm{KL}(\Pi^*)$.*

**Remark 3.** Tokdar [29] assumed that the weak support of $\Pi$ includes all compactly supported probabilities in $\mathbb{R}^d \times \mathbb{R}^+$. Then automatically the weak support of $\Pi$ is $\mathscr{M}(\mathbb{R}^d \times \mathbb{R}^+)$. This is because any arbitrary probability measure can be weakly approximated by a sequence of compactly supported probability measures.

PROOF OF THEOREM 2. We prove this theorem by verifying the conditions in Theorem 1. Since there is no hyper-parameter in Prior 1, we only need to show that Conditions A1 and A3 are met.

To show that Condition A1 is met, we define,

$$f_m(x) = \begin{cases} t_m f_0(x), & \|x\| < m, \\ 0, & \text{otherwise,} \end{cases} \qquad m \geq 1,$$

where $t_m^{-1} = \int_{\|x\| < m} f_0(x) dx$, $h_m = m^{-\eta}$, $F_m$ is the probability measure corresponding to $f_m$, $P_m = F_m \times \delta(h_m)$, where $\delta(\cdot)$ is the degenerate distribution. Obviously, $P_m$ is compactly supported. Then, using the transformation $a = (x - \theta)/h_m$,

$$
\begin{aligned}
f_{P_m}(x) = \int \frac{1}{h_m^d} \chi\left(\frac{x-\theta}{h_m}\right) dF_m(\theta) &= t_m \int_{\|\theta\| < m} \frac{1}{h_m^d} \chi\left(\frac{x-\theta}{h_m}\right) f_0(\theta) d\theta \\
&= \int_{\|x - a h_m\| < m} \chi(a) f_0(x - a h_m) da.
\end{aligned}
$$

Since for any given $a$, $\chi(a)f_0(x - ah_m) \rightarrow \chi(a)f_0(x)$ as $h_m \rightarrow 0$ and $f_0$ is bounded, by the DCT, we obtain $f_{P_m}(x) \rightarrow f_0(x)$.

Now, to satisfy Condition A1, we show that

$$\int f_0(x) \log \frac{f_0(x)}{f_{P_m}(x)} dx \rightarrow 0 \quad \text{as } m \rightarrow \infty.$$

To this end, observe that

$$
\begin{aligned}
f_{P_m}(x) &= t_m \int_{\|\theta\| < m} \frac{1}{h_m^d} \chi \left( \frac{x - \theta}{h_m} \right) f_0(\theta) d\theta \\
&\leq M t_m \int_{\|\theta\| < m} \frac{1}{h_m^d} \chi \left( \frac{x - \theta}{h_m} \right) d\theta \\
&\leq M t_m \leq M t_1.
\end{aligned}
$$

Hence, as $\log \frac{f_0(x)}{M t_1} < 0$,

$$\log \frac{f_0(x)}{f_{P_m}(x)} \geq \log \frac{f_0(x)}{M t_1}. \tag{7}$$

Also

$$
\begin{aligned}
&\int f_0(x) \log \frac{f_0(x)}{f_{P_m}(x)} dx \\
&= \int_{\|x\| \leq m} f_0(x) \log \frac{f_0(x)}{f_{P_m}(x)} dx + \int_{\|x\| > m} f_0(x) \log \frac{f_0(x)}{f_{P_m}(x)} dx.
\end{aligned}
$$

Let $m < l_1$. Now, for $\|x\| > m$, using assumption B2,

$$
\begin{aligned}
f_{P_m}(x) &= t_m \int_{\|\theta\| < m} \frac{1}{h_m^d} \chi \left( \frac{x - \theta}{h_m} \right) f_0(\theta) d\theta \\
&\geq t_m \int_{\|\theta\| < m} \frac{1}{h_m^d} \chi \left( \frac{x + m \frac{x}{\|x\|}}{h_m} \right) f_0(\theta) d\theta \\
&= \frac{1}{h_m^d} \chi \left( \frac{x + m \frac{x}{\|x\|}}{h_m} \right) t_m \int_{\|\theta\| < m} f_0(\theta) d\theta \\
&= \frac{1}{h_m^d} \chi \left( \frac{x + m \frac{x}{\|x\|}}{h_m} \right) \\
&= m^\eta \chi \left( m^\eta x + \frac{x}{\|x\|} m^{1+\eta} \right) \\
&\geq \|x\|^\eta \chi(2\|x\|^\eta x) \tag{8}
\end{aligned}
$$

The last inequality holds when $T \mapsto T^\eta \chi(T^\eta(x + T \frac{x}{\|x\|}))$ is decreasing for $T > T_0$.

This follows because, with $z = T^\eta x + T^{\eta+1} x/\|x\|$, a positive multiple of $x$,

$$\frac{d}{dT}\left\{\eta \log T + \log \chi\left(T^\eta x + T^{\eta+1}\frac{x}{\|x\|}\right)\right\}$$

$$= \frac{\eta}{T} + \sum_{i=1}^{d}\frac{\chi_i'(z)}{\chi(z)}\left(\frac{\eta}{T}z_i + T^\eta\frac{z_i}{\|z\|}\right)$$

$$= \frac{\eta}{T}\left\{1 + \sum_{i=1}^{d}\frac{\chi_i'(z)}{\chi(z)}z_i\left(1 + \frac{T^{1+\eta}}{\eta\|z\|}\right)\right\} \leq 0$$

by Condition B3.

For $\|x\| \leq m$, let $\delta > 0$ be fixed, and $\phi_m^*(x) = \inf_{\|t-x\| < \delta h_m} f_0(t)$,

$$\begin{aligned}
f_{P_m}(x) &= t_m \int_{\|\theta\| \leq m}\frac{1}{h_m^d}\chi\left(\frac{x-\theta}{h_m}\right)f_0(\theta)d\theta \\
&\geq t_m \int_{\{\|\theta\| < m\} \cap \{\|\theta-x\| < \delta h_m\}}\frac{1}{h_m^d}\chi\left(\frac{x-\theta}{h_m}\right)f_0(\theta)d\theta \\
&\geq t_m \phi_m^*(x)\int_{\{\|\theta\| < m\} \cap \{\|\theta-x\| < \delta h_m\}}\frac{1}{h_m^d}\chi\left(\frac{x-\theta}{h_m}\right)d\theta \\
&= t_m \phi_m^*(x)\int_{\{\|x-uh_m\| \leq m\} \cap \{\|u\| \leq \delta\}}\chi(u)du \\
&\geq t_m \phi_m^*(x)\int_{\prod_{i=1}^{d}[0,\text{sign}(x_i)\delta/\sqrt{d}]}\chi(u)du,
\end{aligned}$$

with the convention that $[a, b] = [b, a]$ if $b < a$. The last inequality holds because when $\|x\| \leq m$,

$$\left\{u : u \in \prod_{i=1}^{d}[0,\text{sign}(x_i)\delta/\sqrt{d}]\right\} \subset \left\{u : \|x/h_m - u\| \leq m/h_m \text{ and } \|u\| \leq \delta\right\}.$$

We have $t_m \geq 1$, $\phi_m^*(x) \geq \phi_1(x)$. Let

$$c = \min_{x \in \{\delta/\sqrt{d},\, -\delta/\sqrt{d}\}^d}\left(\int_{\prod_{i=1}^{d}[0,x_i]}\chi(u)du\right).$$

Then, $f_{P_m}(x) \geq c\phi_1(x)$, for all $\|x\| < m$. For $0 < R < m$,

$$f_{P_m}(x) \geq \begin{cases} c\phi_1(x), & \|x\| < R, \\ \min\left\{\|x\|^\eta\chi\left(2\|x\|^{1+\eta}\frac{x}{\|x\|}\right), c\phi_1(x)\right\}, & \|x\| \geq R. \end{cases}$$

$$\log\frac{f_0(x)}{f_{P_m}(x)} \leq \xi(x) := \begin{cases} \log\dfrac{f_0(x)}{c\phi_1(x)}, & \|x\| < R, \\ \max\left(\log\dfrac{f_0(x)}{\|x\|^\eta\chi(2\|x\|^{1+\eta}\frac{x}{\|x\|})}, \log\dfrac{f_0(x)}{c\phi_1(x)}\right), & \|x\| \geq R. \end{cases}$$

$$(9)$$

Combining (7) and (9), we obtain

$$\left| \log \frac{f_0(x)}{f_{P_m}(x)} \right| \leq \max \left( \xi(x), \left| \log \frac{f_0(x)}{Mt_1} \right| \right).$$

From Condition B5,

$$\int \left| \log \frac{f_0(x)}{Mt_1} \right| f_0(x) dx = \log Mt_1 - \int f_0(x) \log f_0(x) dx < \infty.$$

Now

$$\begin{aligned}
\int \xi(x) f_0(x) dx &= \int_{\|x\| < R} f_0(x) \log \frac{f_0(x)}{c\phi_1(x)} dx \\
&\quad + \int_{\|x\| \geq R} f_0(x) \max \left( \log \frac{f_0(x)}{\|x\|^\eta \chi(2\|x\|^\eta x)}, \log \frac{f_0(x)}{c\phi_1(x)} \right).
\end{aligned}$$

Hence,

$$\begin{aligned}
\int \xi(x) f_0(x) dx &\leq \int f_0(x) \log \frac{f_0(x)}{c\phi_1(x)} dx \\
&\quad + \int_{\|x\| \geq R, f_0(x) > \|x\|^\eta \chi(2\|x\|^\eta x)} f_0(x) \log \frac{f_0(x)}{\|x\|^\eta \chi(2\|x\|^\eta x)} dx,
\end{aligned}$$

since $\max(x_1, x_2) \leq x_1 + x_2^+$ if $x_1 \geq 0$. The first term on the RHS of the above inequality is finite, by Condition B6. By Conditions B5 and B7, the second term is also finite. Thus $\int f_0(x) \log \frac{f_0(x)}{f_{P_m}(x)} dx \to 0$ as $m \to \infty$, i.e., Condition A1 is satisfied.

We show that Condition A3 is met by verifying the conditions of Lemma 3. First, from the proof above, we see that for any $\epsilon > 0$, there exists $m_\epsilon$ such that $\int f_0(x) \log \frac{f_0(x)}{f_{P_{m_\epsilon}}(x)} dx < \epsilon$. Let $P_\epsilon$ in Theorem 1 be chosen to be $P_{m_\epsilon}$, which is compactly supported. By Condition B8, $P_\epsilon \in \mathrm{supp}(\Pi)$. Second, Condition A7 is satisfied. To show $\log \frac{f_{P_\epsilon}(x)}{\inf_{(\theta,h) \in D} \frac{1}{h^d} \chi(\frac{x-\theta}{h})}$ is $f_0$-integrable, it suffices to show that $\log f_{P_\epsilon}(x)$ and $\log \inf_{(\theta,h) \in D} \frac{1}{h^d} \chi(\frac{x-\theta}{h})$ are both $f_0$-integrable. Without loss of generality, let $D = \{\|\theta\| \leq a^*\} \times [\underline{h}, \overline{h}]$, where $a^* \geq m_\epsilon$ and $0 < \underline{h} \leq m_\epsilon^{-\eta} \leq \overline{h} < \frac{1}{2}$. For $\|x\| < a^*$, $\log \inf_{(\theta,h) \in D} \frac{1}{h^d} \chi(\frac{x-\theta}{h})$ is bounded. For $\|x\| > a^*$,

$$\log \inf_{(\theta,h) \in D} \frac{1}{h^d} \chi \left( \frac{x-\theta}{h} \right) = \log \left\{ \frac{1}{\underline{h}^d} \chi \left( \frac{x + a^* \frac{x}{\|x\|}}{\underline{h}} \right) \right\}. \tag{10}$$

By Condition B7 and expression (10), $\log \inf_{(\theta,h) \in D} \frac{1}{h^d} \chi(\frac{x-\theta}{h})$ is $f_0$-integrable.

Consider $f_{P_\epsilon}(x) = \int_D \frac{1}{h^d} \chi(\frac{x-\theta}{h}) dP_\epsilon$. Let $D = \{\|\theta\| < a^*\} \times [\underline{h}, \overline{h}]$, then

$$\left| \log \int_D \frac{1}{h^d} \chi \left( \frac{x-\theta}{h} \right) dP_\epsilon \right| \leq \left| \log \left\{ \frac{1}{\underline{h}^d} \chi \left( \frac{x + a^* \frac{x}{\|x\|}}{\underline{h}} \right) P_\epsilon(D) \right\} \right|,$$

for $\|x\| > a^*$. Hence, $\log f_{P_\epsilon}(x)$ is also $f_0$-integrable by the similar argument.

Condition A8 is satisfied by Condition B1.

We show that Condition A9 is also satisfied. Let $C \subset \mathfrak{X}$ be a given compact set. First we show that $\{\frac{1}{h^d}\chi(\frac{x-\theta}{h}) : x \in C\}$ is uniformly equicontinuous as a family of functions of $(\theta, h)$ on $E = [-a, a]^d \times [\frac{1}{2}\underline{h}, 2\overline{h}]$ where $a > a^*$.

Such an $E$ contains $D$ in its interior, and is compact. By the definition of uniform equicontinuity, it is to show that for any $\epsilon > 0$, there exists $\delta > 0$ such that for all $x \in C$ and all $(\theta, h), (\theta', h') \in E$ with $\|(\theta, h) - (\theta', h')\| < \delta$, we have $|h^{-d}\chi(\frac{x-\theta}{h}) - h'^{-d}\chi(\frac{x-\theta'}{h'})| < \epsilon$. Observe that

$$\left| \frac{1}{h^d}\chi\left(\frac{x-\theta}{h}\right) - \frac{1}{h'^d}\chi\left(\frac{x-\theta'}{h'}\right) \right|$$
$$= \left| \frac{\chi(\frac{x-\theta}{h}) - \frac{h^d}{h'^d}\chi(\frac{x-\theta}{h})}{h^d} \right|$$
$$\leq \frac{|\chi(\frac{x-\theta}{h}) - \chi(\frac{x-\theta'}{h'})|}{h^d} + \frac{|h'^d - h^d|}{h^d h'^d}\chi\left(\frac{x-\theta'}{h'}\right). \tag{11}$$

Since $E$ and $C$ are compact and $h$ is bounded away from 0 within $E$, $\{\frac{x-\theta}{h} : x \in C, (\theta, h) \in E\}$ is also a compact set. Hence $c_1 = \sup_{x \in C, (\theta, h) \in E} \chi(\frac{x-\theta'}{h'})$ is finite, by the continuity of $\chi(\cdot)$. Let $\delta^* = \frac{h^{2d}}{2^{2d+1}c_1}\epsilon$, then for $|h' - h| < \frac{\delta^*}{d(2h)^{d-1}}$, we have $|h'^d - h^d| < \delta^*$ and hence the last term in (11) is less than $\epsilon/2$. Since $\{\frac{x-\theta}{h} : x \in C, (\theta, h) \in E$ is compact, $\chi(\cdot)$ is uniformly continuous on it. For any given $\epsilon > 0$, there exists $\delta^{**} > 0$ such that whenever $x \in C$ and $(\theta, h), (\theta', h') \in E$, with $\|\frac{x-\theta}{h} - \frac{x-\theta'}{h'}\| < \delta^{**}$, we have $|\chi(\frac{x-\theta}{h}) - \chi(\frac{x-\theta'}{h'})| < \epsilon \underline{h}/2^{d+1}$, which ensures the second term on the RHS of (11) less than $\epsilon/2$. Notice that $\|\frac{x-\theta}{h} - \frac{x-\theta'}{h'}\| < \delta^{**}$ is equivalent to

$$\|(h - h')\theta + (\theta' - \theta)h + (h' - h)x\| < hh'\delta^{**}. \tag{12}$$

When $\|\theta - \theta'\| < \frac{h^{2d}\delta^{**}}{4h}$ and $|h - h'| < \min\{\frac{h^{2d}\delta^{**}}{4\sqrt{d}a}, |h - h'| < \frac{h^{2d}\delta^{**}}{2\sup_{x \in C}\|x\|}\}$, relation (12) holds. Hence if $\epsilon > 0$ and $\delta = \min\{\frac{\delta^*}{d(2h)^{d-1}}, \frac{h^{2d}\delta^{**}}{24h}, \frac{h^{2d}\delta^{**}}{12\sqrt{d}a}, \frac{h^{2d}\delta^{**}}{12\sup_{x \in C}\|x\|}\}$, then for all $x \in C$ and all $(\theta, h), (\theta', h') \in E$ with $\|(\theta, h) - (\theta', h')\| < \delta$, we have $|h^{-d}\chi(\frac{x-\theta}{h}) - h'^{-d}\chi(\frac{x-\theta'}{h'})| < \epsilon$. Thus the uniform equicontinuity required in Condition A9 is satisfied.

We can enlarge $E$ to ensure that $h^{-d}\chi(\frac{x-\theta}{h})$ is less than any preassigned number for $x \in C$ and $(\theta, h) \in E^c$. This holds for large value of $h$, since $\chi(\cdot)$ is bounded. For small values of $h$, notice that $h^{-d}\chi(\frac{x-\theta}{h}) \leq h^{-d}o(\frac{h^d}{\|x-\theta\|^d}) = o(\|x - \theta\|^{-d})$. This follows from Assumption B9 when $d \geq 2$. For $d = 1$, the condition automatically holds since $\int \chi(y)dy = 1$ implies $\chi(y) = o(\|y\|^{-1})$ with the help of the montonicity condition B2. For given $C$, choosing $a$ and $\overline{h}$ large enough to construct the set $E$, we have $\sup\{h^{-d}\chi(\frac{x-\theta}{h}) : x \in C, (\theta, h) \in E^c\} < c\epsilon/4$, for any given $\epsilon$. $\qquad \square$

Now we consider Prior 2 with location scale family kernels. Let the location-parameter for the density be mixed according to $P$ following a prior $\Pi$. Let the scale-parameter $h$ be a hyper-parameter, which is also given a prior distribution $\mu$. Assume that $h$ and $P$ are a priori independently distributed. We let $\Pi^*$ to denote the prior for the density functions on $\mathfrak{X}$, induced by $\Pi \times \mu$ via the mapping $(P, h) \mapsto f_{P,h} = \int h^{-d} \chi(\frac{x-\theta}{h}) dP(\theta)$. We then have the following theorem.

**Theorem 3.** *For such prior described above, let $\chi(x)$ and $f_0(x)$ be densities on $\mathfrak{X}$ satisfying condition B1–B9. Then, $f_0 \in \mathrm{KL}(\Pi^*)$.*

PROOF. The proof uses Theorem 1 and Lemmas 2 and 3. Verification the Conditions A7–A9 is similar to (but easier than) that in Theorem 2. The second inequality in Condition B7 implies that Condition A5 is satisfied. Conditions A4 and A6 are satisfied since $\chi(\cdot)$ is a continuous probability density function and the kernel we consider here is a location family of $\chi(\cdot)$ with a fixed scale. Condition A1 will be proved in the same way as in the proof of Theorem 2. □

## 4. Examples

In this section, we discuss the KL property for some kernel mixture priors with concretely specified kernels. More precisely, we prove that the property holds under some conditions on the true density when the kernel is chosen to be skew-normal (normal also, as it is a special case), multivariate normal, logistic, double exponential, $t$ (Cauchy also as it is a special case), histogram, triangular, uniform, scaled uniform, exponential, log-normal, gamma, inverse gamma and Weibull densities.

### 4.1. Location-scale kernels

For a given density $\chi(\cdot)$ supported on the entire real line (or $\mathbb{R}^d$ when $\mathfrak{X}$ is $d$-dimensional), we shall consider two types of kernel mixture prior — Prior 1 where both the location parameter $\theta$ and the scale parameter $\phi$ of $\phi^{-d} \chi((x - \theta)/\phi)$ are mixed according to a random probability measure on $\mathbb{R}^d \times (0, \infty)$, or Prior 2 where $\theta$ is mixed according to a random probability measure $P$ on $\mathbb{R}^d$ and $\phi$ is given a prior $\mu$ on $(0, \infty)$. The KL property may be verified by checking Condition B1–B9 for the kernel and applying respectively Theorem 2 or Theorem 3.

In this subsection, we consider several examples of location-scale kernels. Condition B1 and B2 can be easily verified. Conditions B4–B6 are also the conditions assumed in all the following theorems for each of the location scale density kernels. By choosing prior on $P$ as described in Remark 1, Condition B8 can be satisfied. In this subsection, only multivariate normal density has a mixing parameter $\theta$ with dimension $d \geq 2$. For this kernel Condition B9 is obviously satisfied. Hence, in the rest of this subsection, for each kernel function and corresponding prior, we only show that conditions B3 and B7 are satisfied.

## 1. Skew-normal density kernel

Consider the skew-normal kernel

$$\chi_\lambda(x) = 2\frac{1}{\sqrt{2\pi}}e^{-x^2/2}\int_{-\infty}^{\lambda x}\frac{1}{\sqrt{2\pi}}e^{-t^2/2}dt,$$

where the skewness parameter $\lambda$ is given. We have the following result.

**Theorem 4.** *Assume that the prior $\Pi$ satisfies B8. Let $f_0(x)$ be a continuous density on $\mathbb{R}$ satisfying conditions B4, B5, B6 and there exists $\eta > 0$ such that $\int_{\mathbb{R}} |x|^{2(1+\eta)} f_0(x)dx < \infty$. Then $f_0 \in \mathrm{KL}(\Pi^*)$.*

PROOF. For Condition B3, we have $\frac{\chi'(z)}{\chi(z)} = -z + \frac{\Phi'(\lambda z)}{\Phi(\lambda z)}$, $\frac{\Phi'(\lambda z)}{\Phi(\lambda z)} \to \infty$ when $z \to -\infty$ by L'Hospital's rule, since $\frac{(e^{-(\lambda z)^2/2}\lambda)'}{(\int_{-\infty}^{\lambda z} e^{-(\lambda t)^2/2}dt)'} = -\lambda z$; and $\frac{\Phi'(\lambda z)}{\Phi(\lambda z)} \to 0$ when $z \to \infty$. Hence Condition B3 is satisfied.

Condition B7 is satisfied, since

$$\left|\int f_0(x)\log\chi(2|x|^\eta x)dx\right| = \left|\int f_0(x)\left(c_1(x) - \frac{(2|x|^{1+\eta})^2}{2}\right)dx\right| < \infty$$

and similarly

$$\int f_0(x)\left|\log\chi\left(\frac{x-a}{b}\right)\right|dx$$
$$= \int f_0(x)\left|c_2(x) - \frac{(x-a)^2}{2b^2}\right|dx < \infty$$

for any $a$ and $b$, where $c_1(x)$ and $c_2(x)$ are bounded functions here. $\quad\square$

**Remark 4.** With $\lambda = 0$, Theorem 4 implies Theorem 3.2 of [29], since the normal density is a special case of the skew-normal.

## 2. Multivariate normal density kernel

Let $\chi(x) = (2\pi)^{-d/2}\prod_{i=1}^d e^{-x_i^2/2}$, where $x = (x_1, \ldots, x_d)$. We have the following result.

**Theorem 5.** *Assume that the prior $\Pi$ satisfies B8. Let $f_0(x)$ be a continuous density on $\mathbb{R}^d$ satisfying Conditions B4, B5, B6 and that $\int \|x\|^{2(1+\eta)} f_0(x)dx < \infty$ for some $\eta > 0$. Then $f_0 \in \mathrm{KL}(\Pi^*)$.*

PROOF. The proof of this theorem is very similar to the proof of Theorem 4, with $\lambda = 0$ and some other minor modifications in all the steps except in verifying Condition B7. Note that for some bounded functions $c_1(x)$ and $c_2(x)$, we have that

$$\left|\int f_0(x)\log\chi(2\|x\|^\eta x)dx\right|$$
$$= \left|\int c_1(x)f_0(x)dx - \int 2f_0(x)\|x\|^{2(1+\eta)}dx\right| < \infty.$$

and similarly

$$\int f_0(x) \left| \log \chi \left( \frac{x-a}{b} \right) \right| dx = \int f_0(x) \left| c_2(x) - \frac{\sum_1^d (x_i - a_i)^2}{2b^2} \right| dx < \infty$$

for any $a$ and $b$.  □

### 3. Double-exponential density kernel

Let $\chi(x) = \frac{1}{2} e^{-|x|}$. We have the following result.

**Theorem 6.** *Assume that the prior* $\Pi$ *satisfies B8. Let* $f_0(x)$ *be a continuous density on* $\mathbb{R}$ *satisfying B4, B5, B6 and* $\int_{\mathbb{R}} |x|^{1+\eta} f_0(x) dx < \infty$ *for some* $\eta > 0$. *Then* $f_0 \in \mathrm{KL}(\Pi^*)$.

PROOF. Condition B3 is satisfied, since $\frac{\chi'(z)}{\chi(z)} = -1$ when $z > 0$, and $\frac{\chi'(z)}{\chi(z)} = 1$ when $z \leq 0$. Condition B7 follows easily from the fact that $|\log \chi(x)|$ is a linear function of $|x|$.  □

### 4. Logistic density kernel

Let the kernel be $\chi(x) = e^{-x}/(1 + e^{-x})^2$. We have the following result.

**Theorem 7.** *Assume that the prior* $\Pi$ *satisfies B8. Let* $f_0(x)$ *be a continuous density on* $\mathbb{R}$ *satisfying B4, B5, B6 and* $\int_{\mathbb{R}} |x|^{1+\eta} f_0(x) dx < \infty$ *for some* $\eta > 0$. *Then* $f_0 \in \mathrm{KL}(\Pi^*)$.

PROOF. Condition B3 is satisfied, since $\frac{\chi'(z)}{\chi(z)} \to -1$ as $z \to \infty$ and $\frac{\chi'(z)}{\chi(z)} \to 1$ as $z \to -\infty$. Condition B7 is easily verified since the tails of $\log \chi(x)$ behave like $|x|$.  □

### 5. $t_\nu$-density kernel

Let the kernel be given by

$$\chi_\nu(x) = \frac{\Gamma(\frac{\nu+1}{2})}{\sqrt{\nu \pi} \Gamma(\frac{\nu}{2})} \frac{1}{(1 + \frac{(x-\theta)^2}{\phi^2 \nu})^{(\nu+1)/2}},$$

where the degrees of freedom $\nu$ is given. Let $\log_+ u = \max(\log u, 0)$. We have the following result.

**Theorem 8.** *Assume that the prior* $\Pi$ *satisfies B8. Let* $f_0(x)$ *be a continuous density on* $\mathbb{R}$ *satisfying B4, B5, B6 and* $\int_{\mathbb{R}} \log_+ |x| f_0(x) dx < \infty$. *Then* $f_0 \in \mathrm{KL}(\Pi^*)$.

PROOF. Condition B3 is satisfied, since $\frac{\chi'(z)}{\chi(z)} = -cz(1 + \frac{z^2}{\nu})^{-1}$, where $c$ is a positive constant.

Condition B7 can be verified by observing the tail of $|\log \chi_\nu(x)|$ has growth like $\log |x|$ as $|x| \to \infty$.  □

**Remark 5.** Since the Cauchy density is the $t$-density with $\nu = 1$, Theorem 8 applies to the Cauchy kernel.

### *4.2. Kernels with bounded support*

The priors with kernels supported on $[0, 1]$ are preferred for estimating densities supported on $[0, 1]$. We study the KL property of such priors using Theorem 1.

The following lemma will be used in the following proofs repeatedly.

**Lemma 4.** *For any density $f_0$ on $[0,1]$ and $\epsilon > 0$, there exist $m > 0$ and $f_1(x) \geq m > 0$, such that $\Pi^*(\mathscr{K}_\epsilon(f_1)) > 0$ implies that $\Pi^*(\mathscr{K}_{2\epsilon + \sqrt{\epsilon}}(f_0)) > 0$.*

PROOF. If $f_0$ is not bounded away from zero, then define

$$f_1(x) = \frac{\max(f_0(x), m)}{\int \max(f_0(u), m)du}.$$

By Lemma 5.1 in [12], we have $\mathcal{K}(f_0; f) \leq (c+1)\log c + [\mathcal{K}(f_1; f) + \sqrt{\mathcal{K}(f_1; f)}]$, where $c = \int \max(f_0(x), m)dx$. Hence, $c \to 1$ as $m \to 0$. For any given $\epsilon > 0$, there exists $m > 0$ such that $(c+1)\log c < \epsilon$. Therefore $\Pi^*(\mathscr{K}_{2\epsilon + \sqrt{\epsilon}}(f_0)) \geq \Pi^*(\mathscr{K}_\epsilon(f_1))$. □

## 6. Histogram density kernel

Let the kernel function be

$$K(x; \theta, m) = \begin{cases} m, & \text{both } x \text{ and } \theta \in ((i-1)/m, i/m), \text{ for some } 1 \leq i \leq m < \infty, \\ 0, & \text{otherwise.} \end{cases}$$

Consider a kernel mixture prior obtained by mixing both $\theta$ and $m$. We have the following result. An analogous result holds when only $\theta$ is mixed and $m$ is given a prior with infinite support.

**Theorem 9.** *If $f_0(x)$ is a continuous density on $[0, 1]$, and the weak support of $\Pi$ contains $\mathscr{M}([0, 1] \times \mathbb{N})$, then $f_0 \in \mathrm{KL}(\Pi^*)$.*

PROOF. By Lemma 4, we only need to show that Conditions A1 and A3 are satisfied for the density $f_0$ that bounded away from zero. For any $\epsilon > 0$, there exist integer $m > 0$ and $\{w_1, w_2, \cdots, w_m\}$, such that $\sum_{i=1}^m w_i = 1$ and

$$\sup_{x \in [0,1]} \left| f_0(x) - \sum_{i=1}^m w_i K\left(x; \frac{i - \frac{1}{2}}{m}, m\right) \right| < \epsilon. \tag{13}$$

To see this, define $w_i = \frac{f_0(\frac{i-1}{m}) + f(\frac{i}{m})}{\sum_{j=1}^m f_0(\frac{j-1}{m}) + f_0(\frac{j}{m})}$. By Riemann integrability of a continuous function, for any $\epsilon_1 > 0$, there exists $M_1 > 0$, such that for $m > M_1$, $|\sum_1^m \frac{f_0(\frac{i-1}{m}) + f_0(\frac{i}{m})}{2m} - 1| < \epsilon_1$. Since $f_0$ is continuous on a compact set, it is uniformly continuous. Hence, for any given $\epsilon_2 > 0$, there exists $M_2 > 0$, such that for $m > M_2$, $\sup |f_0(x) - \sum_1^m \frac{f_0(\frac{i-1}{m}) + f_0(\frac{i}{m})}{2m} K(x; \frac{i-1/2}{m}, m)| < \epsilon_2$. Let $\Delta = \sum_{i=1}^m \frac{f_0(\frac{i-1}{m}) + f_0(\frac{i}{m})}{2m}$, we have

$$\left| f_0(x) - \sum_{i=1}^m w_i K\left(x; \frac{i - \frac{1}{2}}{m}, m\right) \right| \leq |(\Delta - 1)f_0(x) + \epsilon_2|\frac{1}{\Delta} \leq 2M\epsilon_1 + 2\epsilon_2,$$

where $M$ is an upper bound for $f_0$ on $[0,1]$. Hence, by choosing $\epsilon_1$ and $\epsilon_2$ small enough, there exists $M_3 = \max(M_1, M_2)$ such that for $m > M_3$, (13) holds. Since we consider $f_0$ bounded away from 0 here, Condition A1 will be satisfied by choosing $m_\epsilon$ large enough and appropriate weights $\{w_1, \ldots, w_{m_\epsilon}\}$.

Let

$$\mathscr{W} = \left\{ P : P\left( \left( \frac{i - \frac{1}{2} - \delta_1}{m_\epsilon}, \frac{i - \frac{1}{2} + \delta_1}{m_\epsilon} \right) \times \{m_\epsilon\} \right) > w_i e^{-\epsilon}, \text{ for } i = 1, \ldots, m_\epsilon \right\},$$

where $0 < \delta_1 < 1/4$ and $\epsilon > 0$. Since $\mathscr{W}$ is not empty and it is an open neighborhood of some distribution that belongs to the support of $\Pi$, $P \in \mathscr{W}$, we have with the index $i$ corresponding to the given $x$, $\frac{f_{P_{m_\epsilon}}}{f_P} < e^\epsilon$, and hence $\int f_0 \log \frac{f_{P_{m_\epsilon}}}{f_P} < \epsilon$ for all $P \in \mathscr{W}$. $\qquad\square$

## 7. Triangular density kernel

Let the kernel function be

$$K(x; m, n) = \begin{cases} \begin{cases} 2n - 2n^2 x, & x \in (0, \frac{1}{n}), \\ 0, & \text{otherwise}, \end{cases} & m = 0, \\[2em] \begin{cases} n^2\left(x - \frac{m}{n}\right) + n, & x \in \left(\frac{m-1}{n}, \frac{m}{n}\right), \\ -n^2\left(x - \frac{m}{n}\right) + n, & x \in \left(\frac{m}{n}, \frac{m+1}{n}\right), \\ 0, & \text{otherwise}, \end{cases} & m = 1, 2, \ldots, n-1, \\[3em] \begin{cases} 2n + 2n^2(x-1), & x \in (0, \frac{1}{n}), \\ 0, & \text{otherwise}, \end{cases} & m = n. \end{cases}$$

Construct a kernel mixture prior by mixing both $m$ and $n$. We have the following result.

**Theorem 10.** *Let $f_0(x)$ be a continuous density on $[0,1]$, and the weak support of $\Pi$ contains $\mathscr{M}([0,1] \times \mathbb{N})$. Then $f_0 \in \mathrm{KL}(\Pi^*)$.*

PROOF. Since the mixing parameters are discrete, defining $w_i = \frac{f_0(i/n)}{\sum_{j=0}^n f_0(j/n)}$ and letting $\mathscr{W} = \{P : P(i/n) > w_i e^{-\epsilon}, \text{ for } i = 1, 2, \ldots, n\}$, we can complete the proof as in Theorem 9. $\qquad\square$

## 8. Bernstein polynomial kernel

In the literature, Bernstein polynomials have been used to estimate densities under both frequentist and Bayesian framework. The motivation of the prior comes from the fact that any bounded function on $[0,1]$ can be approximated by a Bernstein polynomial at each point of continuity of the function; see [22].

As in [23; 24], consider a prior $\Pi^*$ induced on $\mathscr{D}(\mathfrak{X})$ by the map

$$(k, (w_0, \ldots, w_k)) \mapsto \sum_{j=0}^k w_j \binom{k}{j} x^j (1-x)^{k-j}$$

and priors $(w_0, \ldots, w_k)|k \sim \Pi_k$ and $k \sim \mu$, where $\mu$ is a discrete distribution supported on the set of all positive integers, $\Pi_k$ is a distribution supported on $(k+1)$-dimensional simplex $\mathscr{P}_k = \{(w_0, \ldots, w_k), \ 0 \le w_j \le 1, \ j = 0, \ldots, k, \ \sum_0^k w_j = 1\}$. We can then rederive Theorem 2 of [26] from Theorem 1.

**Theorem 11.** *If $f_0(x)$ is a continuous density on $[0, 1]$, $\mu(k) > 0$ for infinitely many $k = 1, 2, \ldots$, and $\Pi_k$ is fully supported on $\mathscr{P}_k$, then $f_0 \in \mathrm{KL}(\Pi^*)$.*

PROOF. Though the prior is slightly different from Prior 2 in that $\Pi_k$ is allowed to depend on $k$, we can still use Theorem 1 by changing $\Pi(\mathscr{W}) > 0$ to $\Pi_k(\mathscr{W}) > 0$ for any given $k$. This follows since $k$ is discrete. By Lemma 4, we may assume that $f_0$ is bounded from below. Since Bernstein polynomials uniformly approximate any continuous density (see, for instance, Theorem 1 of [5]), it follows that Condition A1 is satisfied. Condition A3 holds by the discreteness of $k$ and the assumed positivity condition of its prior. The rest of the proof proceeds as before by considering all possible weights $w'_j > w_j e^{-\epsilon}$. □

### *4.3. Kernels supported on $[0, \infty)$*

### 9. Lognormal density kernel

Let the kernel function be $K(x; \theta, \phi) = \frac{1}{\sqrt{2\pi}x\phi} \frac{e^{-(\log x - \theta)^2/(2\phi)^2}}{x}$. Consider a type I or type II mixture prior based on this kernel.

Transform $x \mapsto e^y$ in the kernel function and in $f_0$. If the model using $e^y K(e^y; \theta, \phi)$ as kernel function possess KL property at $e^y f_0(e^y)$, then the corresponding model using $K(x; \theta, \phi)$ as kernel function possess the KL property at $f_0(x)$. This is because of

$$\int_0^\infty f_0(x) \log \frac{f_0(x)}{\int K(x; \theta, \phi) dP(\theta)} dx$$
$$= \int_{-\infty}^\infty e^y f_0(e^y) \log \frac{e^y f_0(e^y)}{e^y \int K(e^y; \theta, \phi) dP(\theta)} dy < \epsilon.$$

For the lognormal kernel, we have the following result.

**Theorem 12.** *Assume that the prior $\Pi$ satisfies B8. Let $f_0(x)$ be a continuous density on $\mathbb{R}^+$ satisfying*

1. *$f_0$ is nowhere zero except at $x = 0$ and bounded above by $M < \infty$;*
2. *$|\int_{\mathbb{R}^+} f_0(x) \log(x f_0(x)) dx| < \infty$;*
3. *$\int_{\mathbb{R}^+} f_0(x) \log \frac{f_0(x)}{\phi_\delta(x)} dx < \infty$ for some $\delta > 0$, where $\phi_\delta(x) = \inf_{|t-x|<\delta} f_0(t)$;*
4. *There exists $\eta > 0$ such that $|\int_{\mathbb{R}^+} f_0(x)| \log x|^{2(1+\eta)} dx| < \infty$.*

*Then $f_0 \in \mathrm{KL}(\Pi^*)$.*

PROOF. Considering the kernel function $\phi^{-1}\chi((y-\theta)/\phi) = \frac{1}{\sqrt{2\pi}\phi} e^{-(y-\theta)^2/(2\phi^2)}$, we can apply Theorem 4 with $\lambda = 0$ or Theorem 5 with $d = 1$. It follows from a change of variable that $g_0(y) := e^y f_0(e^y)$ satisfies B4, B5, B6 and $\int |y|^{2(1+\eta)} g_0(y) dy < \infty$ for some $\eta > 0$. □

## 10. Weibull density kernel

Weibull is a widely used kernel function. Ghosh and Ghosal [17] discussed a model using this density as kernel function and showed posterior consistency useful in survival analysis. However, the assumption for the true density $f_0$ assumed there was quite strong. Here we establish the KL property with this kernel under very general assumptions.

The Weibull kernel is given by $K(x; \theta, \phi) = \theta \phi^{-1} x^{\theta-1} e^{-x^\theta/\phi}$. We can transform this kernel using the map $x = e^y$ to

$$\theta W((y - \theta^{-1} \log \phi)/\theta^{-1}) = e^{\frac{y - \theta^{-1} \log \phi}{\theta^{-1}}} e^{-e^{\frac{y - \theta^{-1} \log \phi}{\theta^{-1}}}},$$

where $W(z) = \exp[z - e^z]$, the location parameter is $\theta^{-1} \log \phi$ and scale parameter is $\theta^{-1}$. We have the following result.

**Theorem 13.** *Let $f_0(x)$ be a continuous density on $\mathbb{R}^+$ satisfying*

1. *$f_0$ is nowhere zero except at $x = 0$ and bounded above by $M < \infty$;*
2. *$|\int_{\mathbb{R}^+} f_0(x) \log(f_0(x)) dx| < \infty$;*
3. *$\int_{\mathbb{R}^+} f_0(x) \log \frac{f_0(x)}{\phi_\delta(x)} dx < \infty$ for some $\delta > 0$, where $\phi_\delta(x) = \inf_{|t-x|<\delta} f_0(t)$;*
4. *there exists $\eta > 0$ such that $e^{2|\log x|^{1+\eta}}$ is $f_0$-integrable;*
5. *the weak support of $\Pi$ contains $\mathscr{M}(\mathbb{R}^+ \times \mathbb{R}^+)$.*

*Then, $f_0 \in \mathrm{KL}(\Pi^*)$.*

PROOF. We need to verify Conditions B3–B7 for kernel $W(\cdot)$ and true density $e^y f_0(e^y)$. Condition B3 is satisfied, since we have $\frac{W'(z)}{W(z)} = 1 - e^z$. To verify Condition B7, observe that Condition 4 of this theorem implies

$$\left| \int_{\mathbb{R}} e^y f_0(e^y) \log e^{2|y|^{1+\eta}} W(e^{2|y|^{1+\eta}}) dy \right| < \infty$$

and

$$\int_{\mathbb{R}} e^y f_0(e^y) |\log W(e^{\frac{y-a}{b}})| dy < \infty.$$

$\square$

## 11. Gamma density kernel

The gamma density is one of the most widely used kernel function for density estimation on $[0, \infty)$. Hason [18] discussed a model using the gamma density as kernel with the hierarchical structure has as many stages as the most general one we discussed in Section 1. Chen [4] and Bouezmarni and Scaillet [3] discussed a mixture of gamma model with a different parametrization.

Let $K(x; \alpha, \beta) = \frac{1}{\Gamma(\alpha)\beta^\alpha} x^{\alpha-1} e^{-x/\beta}$ be the kernel function. Set

$$\phi_\delta(x) = \begin{cases} \inf_{[x, x+\delta)} f_0(t), & 0 < x < 1, \\ \inf_{(x-\delta, x]} f_0(t), & x \geq 1. \end{cases} \tag{14}$$

**Theorem 14.** *Assume that the weak support of prior $\Pi$ is $\mathscr{M}(\mathbb{R}^+ \times \mathbb{R}^+)$. Let $f_0(x)$ be a continuous and bounded density on $[0, \infty)$ satisfying B4, B5 and*

*B6.\** $\int f_0(x) \log \frac{f_0(x)}{\phi_\delta(x)} dx < \infty$ *for some $\delta > 0$;*
*B7.\* there exists $\eta > 0$, such that $\int \max(x^{-\eta-2}, x^{\eta+2}) f_0(x) dx < \infty$.*

*Then, $f_0 \in \mathrm{KL}(\Pi^*)$.*

PROOF. We use $K_m(x; \alpha)$ to denote $K(x; \alpha, m^{-1})$. Let

$$f_m(x) = t_m \int_2^{1+m^2} K_m(x; \alpha) m^{-1} f_0((\alpha-1)/m) d\alpha, \tag{15}$$

where $t_m = (\int_{m^{-1}}^m f_0(s) ds)^{-1}$. Let $P_m$ denote $F_m^* \times \delta(m^{-1})$, where $F_m^*$ is the probability measure corresponding to $t_m m^{-1} f_0((\alpha-1)/m) \mathbb{1}(\alpha \in [2, 1+m^2])$ as a density function for $\alpha$, and $\mathbb{1}(\cdot)$ is the indicator function. Obviously, $P_m$ is compactly supported and $f_m(x) = f_{P_m}(x)$. Let $F_m$ be the probability measure corresponding to $f_m$. By Lemma 5 in the Appendix, $\int f_0(x) \log \frac{f_0(x)}{f_m(x)} dx \to 0$ as $m \to \infty$, which implies that Condition A1 is satisfied.

To complete the proof, we show that Condition A3 is satisfied by verifying conditions of Lemma 3. For any given $\epsilon > 0$, let $D = [2, 1+m_\epsilon^2] \times \{m_\epsilon^{-1}\}$, where $m_\epsilon$ is such that $\int f_0(x) \log \frac{f_0(x)}{f_{m_\epsilon}(x)} dx < \epsilon$. To verify Condition A7, it is sufficient to show that $\int f_0(x) |\log f_{m_\epsilon}(x)| dx < \infty$ and $\int f_0(x) |\log \inf_{(\alpha,\beta) \in D} K(x; \alpha, \beta)| dx < \infty$. Based on expression (19), (20) and (25) in the appendix, we have

$$\log \inf_{(\alpha,\beta) \in D} K(x; \alpha, \beta) = \log(\min\{K(x; 1+m_\epsilon^2, m_\epsilon^{-1}), K(x; 2, m_\epsilon^{-1})\}),$$

for any $0 < x < \infty$. Hence

$$\left| \log \inf_{(\alpha,\beta) \in D} K(x; \alpha, \beta) \right|$$
$$< xm_\epsilon + (m_\epsilon^2)|\log x| + \left| \log \left(\Gamma(m_\epsilon^2 + 1) m_\epsilon^{-(m_\epsilon^2+1)}\right) \right| + |\log(m_\epsilon^{-2})|.$$

By Condition B7\*, we have that $\int |\log \inf_{(\alpha,\beta) \in D} K(x; \alpha, \beta)| f_0(x) dx < \infty$. Further, $\log f_{m_\epsilon}(x)$ is also $f_0$-integrable by a similar argument. Condition A8 is obviously satisfied. Condition A9 is satisfied by letting $E$ be large enough compact set containing $D$. This proves the theorem. □

## 12. Inverse gamma density kernel

The inverse gamma density function is defined as $h(x; a, b) = \frac{b^a}{\Gamma(a)} x^{-a-1} e^{-b/x}$. We consider the following reparametrerization $K(x; k, z) = h(x; k, kz)$ as the kernel function and construct mixture priors. Let $\phi_\delta$ defined as in (14). We have the following result.

**Theorem 15.** *Assume that the weak support of prior $\Pi$ contains $\mathscr{M}(\mathbb{R}^+ \times \mathbb{R}^+)$. Let $f_0(x)$ be a continuous and bounded density on $[0, \infty)$ satisfying B4, B5, B6\* and B7\*. Then $f_0 \in \mathrm{KL}(\Pi^*)$.*

PROOF. Observe that

$$
\begin{aligned}
\int_0^\infty h(x;k,kz)dP(z) &= \int_0^\infty \frac{(kz)^k}{\Gamma(k)} x^{-(k+1)} e^{-kz/x} dP(z) \\
&= \int_0^\infty \frac{(k/x)^{k+1}}{\Gamma(k+1)} z^k e^{-(k/x)z} dP(z) \\
&= \int_0^\infty g(z;k+1,x/k)dP(z),
\end{aligned}
$$

where $g$ is the gamma density. By Proposition 3.1 in [3], we have for any $x \in [0,\infty)$, $\int_0^\infty g(z;k+1,x/k)f_0(z)dz \to f_0(x)$ as $k \to \infty$, i.e., $\int_0^\infty K(x;k,z)f_0(z)dz \to f_0(x)$ as $k \to \infty$.

Set $f_m(x) = t_m \int_{m^{-1}}^m K(x;k,z)f_0(z)dz$, where $t_m = (\int_{m^{-1}}^m f_0(z)dz)^{-1}$, and let $P_m$ be $F_m \times \delta(m)$, where $F_m$ is the probability measure corresponding to $t_m f_0(x) \mathbb{1}(x \in [m^{-1}, m])$.

Observe that $\frac{d}{dz} \log(h(x;m,mz)) = m(z^{-1} - x^{-1})$. Hence

$$
h(x;m,mz) \geq
\begin{cases}
\dfrac{1}{\Gamma(m)} x^{-m-1} e^{x^{-1}}, & \text{for } x > m, \\[2mm]
\dfrac{m^{2m}}{\Gamma(m)} x^{-m-1} e^{-m^2/x}, & \text{for } x < m^{-1}.
\end{cases}
$$

The derivative of the logarithm of the expression on the RHS of above relation are given by,

$$
\frac{d}{dm} \log\left( \frac{x^{-m-1} e^{-1/x}}{\Gamma(m)} \right) = -\log x - \Psi_0(m) < 0,
$$

for $x > m$, and

$$
\frac{d}{dm} \log\left( \frac{m^{2m} x^{-m-1} e^{-m^2/x}}{\Gamma(m)} \right) = 2\log m + 2 - \log x - \frac{2m}{x} - \Psi_0(m) < 0,
$$

for $x < m^{-1}$, where $\Psi_0(\cdot)$ is the digamma function, and its details is given in the proof of Lemma 5 in the Appendix. Therefore

$$
h(x;m,mz) \geq
\begin{cases}
\dfrac{1}{\Gamma(x)} x^{-x-1} e^{-x^{-1}}, & \text{for } x > m, \\[2mm]
\dfrac{x^{-2/x}}{\Gamma(x^{-1})} x^{-x^{-1}-1} e^{-x^{-3}}, & \text{for } x < m^{-1}.
\end{cases}
$$

and hence

$$
f_m(x) \geq
\begin{cases}
\dfrac{1}{\Gamma(x)} x^{-x-1} e^{-x^{-1}}, & \text{for } x > m, \\[2mm]
\dfrac{x^{-2/x}}{\Gamma(x^{-1})} x^{-x^{-1}-1} e^{-x^{-3}}, & \text{for } x < m^{-1}.
\end{cases}
$$

Now, for $m^{-1} < x < 1$

$$\int_{x \vee m^{-1}}^{x+\delta} \frac{m^m z^m}{\Gamma(m) x^{m+1}} e^{-mz/x} dz = \frac{1}{\Gamma(m+1)} \int_{m \vee x^{-1}}^{m(1+\delta/x)} \left(\frac{mz}{x}\right)^m e^{-\frac{mz}{x}} d\frac{mz}{x}$$

$$\geq \mathrm{Ga}\left((\delta/x+1)m\right) - \mathrm{Ga}(m), \qquad (16)$$

where Ga is the cumulative distribution function (c.d.f.) of gamma distribution with parameter $(m+1, 1)$. For large $m$, the last expression is bounded below by $\{\Phi(1+\delta/x) - \Phi(1)\}/2$ in view of the central limit theorem. Similarly, for $1 \leq x < m$ and large $m$, the lower bound for the left hand side (LHS) of (16) is $\{\Phi(1) - \Phi(1-\delta/x)\}/2$.

Let

$$C(x) = \begin{cases} \{\Phi(1+\delta/x) - \Phi(1)\}/2, & 0 < x < 1, \\ \{\Phi(1) - \Phi(1-\delta/x)\}/2, & x \geq 1. \end{cases}$$

Now we have that $f_m(x) \geq \phi_\delta(x) C(x)$ and $\int |\log C(x)| f_0(x) dx < \infty$. As in the proof of Lemma 5 in the Appendix, this gives a lower bound of $f_m(x)$ for $m^{-1} \leq x \leq m$.

Now we have the lower bound of $f_m(x)$ for any large $m$. Along the same line as for gamma kernel in Lemma 5, we can show that $\int f_0(x) \log \frac{f_0(x)}{f_m(x)} dx \to 0$ as $m \to \infty$, which implies that Condition A1 is satisfied. Similarly as in the proof of Theorem 14, we can show that Condition A3 is also satisfied. $\qquad \square$

## 13. Exponential density kernel

Consider a mixture prior based on the exponential kernel. Let $K(x; \theta) = \theta e^{-\theta x}$. Recall that a function $\varphi$ on $\mathbb{R}^+$ is completely monotone if it possesses derivatives $\varphi^{(n)}$ of all orders and $(-1)^n \varphi^{(n)}(x) \geq 0$ for $x > 0$. Let $\bar{F}_0(x) = 1 - F_0(x)$, where $F_0$ is the distribution function corresponding to density function $f_0$. We have the following result.

**Theorem 16.** *If $f_0$ is a continuous density on $\mathbb{R}^+$, $x$ and $|\log f_0(x)|$ are $f_0$-integrable, $\bar{F}_0(x)$ is completely monotone, and the weak support of $\Pi$ is $\mathscr{M}(\mathbb{R}^+)$, then $f_0 \in \mathrm{KL}(\Pi^*)$.*

PROOF. Since $\bar{F}_0(x)$ is completely monotone, by Theorem 1 in [8, Chapter XIII.4], it is the Laplace transform of a probability distribution $P_0$, i.e., $\bar{F}_0(x) = \int_0^\infty e^{-\theta x} dP_0(\theta)$. Taking derivative on both sides,

$$f_0(x) = -\frac{d}{dx} \int_0^\infty e^{-\theta x} dP_0(\theta) = \int_0^\infty \theta e^{-\theta x} dP_0(\theta) = \int_0^\infty K(x; \theta) dP_0(\theta).$$

Hence under the conditions in this theorem, the true density is of the form of mixture of the kernel.

Let $P_a(A) = P_0(A \cap [a^{-1}, a]) / P_0([a^{-1}, a])$ for any $A \subset \mathbb{R}^+$ and $f_{P_a}$ denote $\int_0^\infty K(x; \theta) dP_a(\theta)$. For any $x \in (0, \infty)$, $\int_0^\infty K(x; \theta) dP_a(\theta) \to \int_0^\infty K(x; \theta) dP_0(\theta)$ as $a \to \infty$. Hence, $\log \frac{f_0(x)}{f_{P_a}(x)} \to 0$ pointwise as $a \to \infty$.

Since $|\log f_0(x)|$ is $f_0$-integrable, showing $|\log f_{P_a}(x)|$ is not greater than an $f_0$-integrable function suffices for an application of DCT to obtain that $\int_0^\infty f_0(x) \log \frac{f_0(x)}{f_{P_a}(x)} dx \to 0$ as $a \to \infty$.

Note that $f_{P_a}(x) = \frac{1}{P_0([a^{-1},a])} \int_{a^{-1}}^a K(x,\theta) dP_0(\theta)$ and $f_0(x) = \int K(x,\theta) dP_0(\theta)$. There exists $a_0 > 0$ such that for $a > a_0$, $f_{P_a}(x) < 2f_0(x)$.

Observe that for given $x$, $K(x,\theta)$ is increasing on $(0, x^{-1}]$ and decreasing on $[x^{-1}, \infty)$ as a function of $\theta$. We obtain the lower bound of $f_{P_a}(x)$ by using this property. First, let $\theta_1$, $\theta_2$ and $\theta_3$ be such that $P_0((0,\theta_1)) = q_1 > 0$, $P_0((\theta_1, \theta_2)) = q_2 > 0$, $P_0((\theta_2, \theta_3)) = q_3 > 0$, $P_0((\theta_3, \infty)) = q_4 > 0$. Choose $a$ sufficiently large, such that $a^{-1} < \theta_1$ and $a > \theta_3$.

For $x \geq \theta_2^{-1}$, $K(x,\theta)$ is decreasing as a function of $\theta$ on $[\theta_2, \infty)$. Hence $f_{P_a}(x) > \theta_3 e^{-x\theta_3} q_3$. For $0 < x < \theta_2^{-1}$, $K(x;\theta)$ is increasing as a function of $\theta$ on $(0, \theta_2)$. Hence, $f_{P_a}(x) > \theta_1 e^{-x\theta_1} q_2$.

Therefore, for $a$ large, we have

$$2f_0(x) > f_{P_a}(x) > \begin{cases} \theta_3 e^{-x\theta_3} q_3, & x \geq 1/\theta_2, \\ \theta_1 e^{-x\theta_1} q_2, & x < 1/\theta_2. \end{cases}$$

Hence,

$$\left| \log \frac{f_0(x)}{f_{P_a}(x)} \right| \leq |\log f_0(x)| + |\log f_{P_a}(x)|,$$

and

$$|\log f_{P_a}(x)| < \max\{\log 2 + |\log f_0(x)|, |\log(\theta_3 q_3)| + |x\theta_3|, |\log(\theta_1 q_2)| + |x\theta_1|\}.$$

Since $\log f_0(x)$ and $x$ are both $f_0$-integrable, by the DCT, we have

$$\int_0^\infty f_0(x) \log \frac{f_0(x)}{f_{P_a}(x)} dx \to 0$$

as $a \to \infty$. Thus Condition A1 is satisfied.

To show that Condition A3 is satisfied, we verify that Conditions A7–A9 are satisfied. For any $\epsilon > 0$, there exists $a > 1$ such that $\int_0^\infty f_0(x) \log \frac{f_0(x)}{f_{P_a}(x)} dx < \epsilon$. From above, we have that $\int \log(f_{P_a}(x)) f_0(x) dx < \infty$. Let $D = [a^{-1}, a]$, then $|\log(\inf_{\theta \in D} K(x;\theta))| \leq xa^{-1} + xa + \log a$. By DCT, $\log(\inf_{\theta \in D} K(x;\theta))$ is $f_0$-integrable. Hence, Condition A7 is satisfied. Condition A8 holds obviously. For Condition A9, the uniform equicontinuity holds for any compact $E$. Without loss of generality, let $C = [c_1, c_2]$, $E = [(ab)^{-1}, ab]$, where $b > 1$, and hence $E \supset D$. Choosing $b$ such that $(ab)^{-1} < c_2^{-1}$ and $ab > c_1^{-1}$, then, by the monotonicity property of exponential density function, $\sup\{K(x,\theta) : x \in C, \theta \in E^c\} = \max(\frac{1}{ab} e^{-\frac{1}{ab}c_1}, abe^{-abc_1}) \to 0$ as $b \to \infty$, so A9 is satisfied. Thus, $f_0 \in \mathrm{KL}(\Pi^*)$. □

## 14. Scaled uniform density kernel

Let the true density $f_0$ be supported on $\mathfrak{X} = \mathbb{R}^+$, and consider a mixture prior based on the scaled uniform kernel $K(x;\theta) = \theta^{-1} \mathrm{1\!l}\{0 \leq x \leq \theta\}$.

**Theorem 17.** *If $f_0(x)$ is a continuous and decreasing density function on $\mathbb{R}^+$ such that $\int f_0 |\log f_0| < \infty$ and the weak support of $\Pi$ is $\mathscr{M}(\mathbb{R}^+)$, then $f_0 \in$ KL$(\Pi^*)$.*

PROOF. We will show that Conditions A1 and A3 are satisfied. Let $x_1 > 0$ and $x_2 > 0$, such that $f_0(x_1) = a$ and $f_0(x_2) = b$, where $0 < b < 1$ and $b < a < f_0(0)$. For given $m$, let $m_1$ and $m_2$ be such that $\frac{m_1}{m} \leq x_1 \leq \frac{m_1+1}{m}$ and $\frac{m_2}{m} \leq x_2 \leq \frac{m_2+1}{m}$.

Let

$$
w_i^* = \begin{cases}
\dfrac{i}{m}\left(f_0(\dfrac{i}{m}) - f_0(\dfrac{i+1}{m})\right), & 1 \leq i < m_1, \\[2mm]
\dfrac{m_1}{m}\left(f_0\left(\dfrac{m_1}{m}\right) - a\right), & i = m_1, \\[2mm]
\dfrac{(m_1+1)}{m}\left(a - f_0\left(\dfrac{m_1+1}{m}\right)\right), & i = m_1 + 1, \\[2mm]
\dfrac{i}{m}\left(f_0\left(\dfrac{i-1}{m}\right) - f_0(\dfrac{i}{m})\right), & m_1 + 1 < i \leq m_2, \\[2mm]
\dfrac{i}{m}\left(f_0\left(\dfrac{i-1}{m}\right) - f_0(\dfrac{i}{m})\right), & i \geq m_2 + 1.
\end{cases}
$$

We define $f_m^*(x) = \sum_1^\infty w_i^* K(x; \frac{i}{m})$. By the continuity of $f_0$, $f_m^*$ converges to $f_0$ pointwise. Note that $f_m^*$ is not a p.d.f. Let $w_i = w_i^* \frac{1-\sum_1^{m_1} w_i - \sum_{m_2+1}^\infty w_i}{\sum_{m_1}^{m_2} w_i^*}$, for $m_1 \leq i \leq m_2$ and $w_i = w_i^*$ for all other $i$'s. Then $w_i$'s are positive and $\sum_1^\infty w_i = 1$. Let $f_m(x) = \sum_1^\infty w_i K(x; \frac{i}{m})$. Observe that

$$
\begin{aligned}
f_m^*(x) - f_m(x) &= \left(\frac{1 - \sum_1^{m_1} w_i - \sum_{m_2+1}^\infty w_i}{\sum_{m_1}^{m_2} w_i^*} - 1\right)\left(\sum_{m_1}^{m_2} w_i^* \frac{m}{i}\right) \\
&\leq \left(\frac{1 - \sum_1^{m_1} w_i - \sum_{m_2+1}^\infty w_i}{\sum_{m_1}^{m_2} w_i^*} - 1\right)\left(\sum_{m_1}^{m_2} w_i^*\right)\frac{m}{m_1} \\
&= \left(1 - \sum_1^{m_1} w_i - \sum_{m_2+1}^\infty w_i - \sum_{m_1}^{m_2} w_i^*\right)\frac{m}{m_1} \\
&= \left(1 - \frac{1}{m}\sum_1^\infty f(i/m) - \frac{a}{m}\right)\frac{m}{m_1} \to 0
\end{aligned}
\tag{17}
$$

as $m \to \infty$, by the definition of Riemann integral. Thus $f_m$ converges to $f_0$ pointwise. Let $m$ large such that the expression on the RHS of (17) is less than $\frac{a}{2}$, we have that

$$
|\log f_m(x)| \leq \begin{cases}
\max(\log 2 + |\log f_0|, |\log a - \log 2|), & 0 < x \leq m_1 + 1, \\
\max(\log a, \log(f_0(x_2 + 1))), & m_1 + 1 < x \leq m_2 + 1, \\
|\log f_0|, & x > m_2 + 1.
\end{cases}
$$

Since $|\log f_0(x)|$ is $f_0$-integrable, by using the DCT, we have $\int f_0 \log \frac{f_0}{f_m} \to 0$ as $m \to \infty$. Condition A3 is satisfied by a similar argument as in the proof for Theorem 9. $\qquad\square$

## Appendix A

**Lemma 5.** *Let $f_m(x)$ be defined as in (15). If the conditions of Theorem 14 are satisfied, then $\int f_0(x) \log \frac{f_0(x)}{f_m(x)} dx \to 0$ as $m \to \infty$.*

PROOF. First, we derive the lower bound of $f_m(x)$ for $x$ in different intervals. Observe that

$$\frac{d}{d\alpha} \log(K_m(x; \alpha)) = \log m + \log x - \Psi_0(\alpha), \tag{18}$$

where $\Psi_0(z) = \frac{d}{dz} \log(\Gamma(z))$, is the digamma function. Also $\Psi_0(z)$ is continuous and monotone increasing for $z \in (0, \infty)$, $\Psi_0(z+1) = \Psi_0(z) + \frac{1}{z}$, and $\Psi_0(z) - \log(z-1) \to 0$; see [2, pp. 549–555] for details.

For $x < m^{-1}$, $\log(mx) < 0$, and $\Psi_0(\alpha) \geq \Psi_0(2) = 0.42$ for $\alpha \in [2, 1+m^2]$, and hence $\frac{d}{d\alpha} \log(K_m(x; \alpha)) < 0$. For $x > m + m^{-1}$ and $\alpha \in [2, 1 + m^2]$, $\log(mx) \geq \log(m^2) \geq \Psi_0(1+m^2) \geq \Psi_0(\alpha)$, and hence $\frac{d}{d\alpha} \log(K_m(x; \alpha)) > 0$. Thus replacing $\alpha$ by $1 + m^2$ in the integrand, we obtain a lower bound for $f_m(x)$, $x < m^{-1}$, as,

$$f_m(x) \geq t_m \int_2^{1+m^2} \frac{x^{m^2} e^{-xm} m^{m^2+1}}{\Gamma(m^2+1)} f_0(\alpha) d\alpha = \frac{x^{m^2} e^{-xm} m^{m^2+1}}{\Gamma(m^2+1)}. \tag{19}$$

Similarly, replacing $\alpha$ by 2 in the integrand, we obtain that for $x > m + m^{-1}$,

$$f_m(x) \geq xe^{-xm}m^2. \tag{20}$$

Consider the RHS of equation (19). For $x < m^{-1}$, we have

$$\frac{d}{dm} \log\left(\frac{x^{m^2} e^{-xm} m^{m^2+1}}{\Gamma(m^2+1)}\right) = 2m[\log(xm) - \Psi_0(m^2+1)] + \frac{m^2+1}{m} - x < 0,$$

for all $m$ sufficiently large, where $c_1 > 0$ is some constant. Consider the RHS of equation (20), for $x > m + m^{-1}$, we have $\frac{d}{dm}\left(xe^{-xm}m^2\right) = xme^{-xm}(2 - xm) < 0$.

Hence, replacing $m$ by $x^{-1}$ on the RHS of (19), we obtain a lower bound of $f_m(x)$ for $x < m^{-1}$ as below,

$$f_m(x) \geq \frac{x^{m^2} e^{-xm} m^{m^2+1}}{\Gamma(m^2+1)} \geq \frac{x^{x^{-2}} e^{-1} x^{-x^{-2}-1}}{\Gamma(x^{-2}+1)} = \frac{1}{ex\Gamma(x^{-2}+1)}; \tag{21}$$

and similarly, replacing $m$ by $x$ on the RHS of (20), we obtain that for $x > m + m^{-1}$,

$$f_m(x) \geq xe^{-xm}m^2 \geq e^{-x^2}x^3. \tag{22}$$

Now, we consider $f_m(x)$ for $m^{-1} \leq x \leq m + m^{-1}$. Let $\delta > 0$ be fixed and $v = (\alpha - 1)/m$. For $m$ large,

$$
\begin{aligned}
f_m(x) &\geq \int_{x-\delta}^{x+\delta} K_m(x; mv+1) t_m f_0(v) dv \\
&\geq \begin{cases} \phi_\delta(x) t_m \displaystyle\int_{m^{-1}\vee x}^{x+\delta} K_m(x; mv+1) dv, & x < 1 \\ \phi_\delta(x) t_m \displaystyle\int_{x-\delta}^{m\wedge x} K_m(x; mv+1) dv, & x \geq 1 \end{cases} \\
&\geq C(x)\phi_\delta(x),
\end{aligned}
$$

where $C(x)$ is given in Lemma 8.

Now we have the lower bound of function $f_m(x)$,

$$
f_m(x) \geq \begin{cases} C(x)\phi_\delta(x), & R^{-1} \leq x \leq R, \\ \min\left(C(x)\phi_\delta(x), \dfrac{1}{ex\Gamma(x^{-2}+1)}\right), & 0 < x < R^{-1}, \\ \min(C(x)\phi_\delta(x), e^{-x^2}x^3), & R < x, \end{cases} \tag{23}
$$

where $0 < R < m$. Hence, we have that

$$
\log \frac{f_0(x)}{f_m(x)} \leq \xi(x)
$$

$$
:= \begin{cases} \log \dfrac{f_0(x)}{C(x)\phi_\delta(x)}, & R^{-1} \leq x \leq R, \\ \max\left\{ \log \dfrac{f_0(x)}{C(x)\phi_\delta(x)}, \log(ex\Gamma(x^{-2}+1)f_0(x)) \right\}, & 0 < x < R^{-1}, \\ \max\left\{ \log \dfrac{f_0(x)}{C(x)\phi_\delta(x)}, \log \dfrac{f_0(x)}{e^{-x^2}x^3} \right\}, & R < x. \end{cases}
$$

Since $f_0(x) < M < \infty$, we also have that $\log \frac{f_0}{f_m} \geq \log \frac{f_0(x)}{Mt_2}$. Further, as $\log \frac{f_0(x)}{Mt_2} < 0$, we have $|\log \frac{f_0(x)}{f_m(x)}| \leq \max\{\xi(x), |\log \frac{f_0(x)}{Mt_2}|\}$.

By Condition B5, $\int |\log \frac{f_0(x)}{Mt_2}| f_0(x) dx = \log Mt_2 - \int f_0 \log(f_0) dx < \infty$. Now, consider $\int \xi(x) f_0(x) dx$, which equals to

$$
\begin{aligned}
&\int_{R^{-1}}^{R} f_0(x) \log \frac{f_0(x)}{C(x)\phi_\delta(x)} dx \\
&+ \int_0^{R^{-1}} f_0(x) \max\left\{ \log \frac{f_0(x)}{C(x)\phi_\delta(x)}, \log(f_0(x)) - \log(ex\Gamma(x^{-2}+1)f_0(x)) \right\} dx \\
&+ \int_R^{\infty} f_0(x) \max\left\{ \log \frac{f_0(x)}{C(x)\phi_\delta(x)}, \log(f_0(x)) - \log \frac{f_0(x)}{e^{-x^2}x^3} \right\} dx \\
&\leq \int_0^{\infty} f_0(x) \log \frac{f_0(x)}{\phi_\delta(x)} dx + \int_0^{\infty} f_0(x) \log \frac{1}{C(x)} dx \\
&\quad + \int_{(0,R^{-1}]\cap A} f_0(x) \left[ \log(ex\Gamma(x^{-2}+1)f_0(x)) \right] dx \\
&\quad + \int_{(R,\infty)\cap B} f_0(x) \left[ \log \frac{f_0(x)}{e^{-x^2}x^3} \right] dx,
\end{aligned} \tag{24}
$$

where $A = \{x : f_0(x) \geq [ex\Gamma(x^{-2} + 1)]^{-1}\}$, and $B = \{x : f_0(x) \geq e^{-x^2}x^3\}$. The above relation (24) holds since $C(x) < 1$ by Lemma 8 and $\max(x_1, x_2) \leq x_1 + x_2^+$ if $x_1 > 0$.

The first term on the RHS of (24) is less than infinity by Condition B6*. By Lemma 8, the second terms on the RHS of (24) is also less than infinity. Note that, by Stirling's inequality, (see [8, vol. I. pp. 50–53])

$$\left| \log \frac{1}{ex\Gamma(x^{-2} + 1)} \right|$$
$$\leq \ |\log x| + 1 + \log(2\pi) + (x^{-2} + 1)\log(x^{-2} + 1) + \frac{(x^{-2} + 1)^2 + 1}{12(x^{-2} + 1)},$$

for $0 < x < 1$. Hence, the third term on the RHS of (24) is less than infinity by Condition B7*. Similarly, so is the fourth term. By Lemma 6, we have that $f_m \to f_0$ pointwise. Thus, by the DCT, $\int f_0(x) \log \frac{f_0(x)}{f_m(x)} dx \to 0$ as $m \to \infty$. $\square$

**Lemma 6.** *Let $f_m(x)$ be defined as in (15), then $f_m(x) \to f_0(x)$ as $m \to \infty$ for each $x > 0$.*

To prove this lemma, we need the lemma below, which generalizes Theorem 2.1. of Devore and Lorentz (1993) from two aspects — the functions $K_m$ and $f$ are considered on a possibly non-compact $\mathfrak{X}$, and the intervals $A_m$ can vary with $m$.

**Lemma 7.** *Let $A_m = [a_m, b_m] \subset \mathfrak{X}$, and let $K_m(x; t)$ be a sequence of continuous functions for $x \in \mathfrak{X}$ and $t \in A_m$. Define $f_m(x) = \int_{A_m} K_m(x, t) f(t) dt$, $m = 1, 2, \ldots$, where $f$ is bounded, uniformly continuous and integrable on $\mathfrak{X}$. If $K_m$ satisfies*

*C1. $\int_{A_m} K_m(x, t) dt \to 1$ as $m \to \infty$,*

*C2. for each $\delta > 0$, $\int_{|x-t| \geq \delta, t \in A_m} |K_m(x, t)| dt \to 0$ as $m \to \infty$,*

*C3. $\int_{A_m} |K_m(x, t)| dt \leq M(x) < \infty$ for each $x \in \mathfrak{X}$, $m = 1, 2 \ldots$, where the bound $M(x)$ may depend on $x$,*

*then $f_m(x) \to f(x)$ for each $x \in \mathfrak{X}$.*

PROOF. Let $\epsilon > 0$ be given and let $\delta > 0$ be so small that $|f(t) - f(x)| < \epsilon$ for $|x - t| \leq \delta$. Because of Condition C1,

$$f_m(x) - f(x) = \int_{A_m} [f(t) - f(x)] K_m(x, t) dt + o(1),$$

where the last term goes to 0 for $m \to \infty$, for each $x \in \mathfrak{X}$. We have

$$\left| \int_{|x-t| \leq \delta, t \in A_m} [f(t) - f(x)] K_m(x, t) dt \right| \leq \epsilon \int_{|x-t| \leq \delta, t \in A_m} |K_m(x, t)| dt \leq \epsilon M(x).$$

It follows from Condition C2 that for each $\delta > 0$, and any bounded continuous function $f^*$ on $\mathfrak{X}$, $\int_{|x-t| \geq \delta} f^*(t) K_m(x, t) dt \to 0$ as $m \to \infty$. Hence,

$$\int_{|x-t| > \delta, t \in A_m} [f(t) - f(x)] K_m(x, t) dt \to 0 \text{ as } m \to \infty.$$

By Condition C3 it now follows that $|f_m(x) - f(x)| \leq \epsilon M(x) + o(1)$, and hence the result. $\square$

PROOF OF LEMMA 6. Let $v = (\alpha - 1)m^{-1}$ and $u = m^{-1}$. Let

$$K(x; v, u) = \frac{x^{v/u}e^{-x/u}}{\Gamma(v/u + 1)u^{v/u+1}},$$

and $K_m(x; v) = K(x; v, m^{-1})$, where $v \in A_m$, $A_m = [m^{-1}, m]$. Now $f_m(x) = \int_{m^{-1}}^{m} K_m(x, v) f_0(v) dv$, we show that such $K_m(x, v)$ satisfies condition C1–C3 in Lemma 7.

Given $x > 0$, consider expression (18), for $m$ sufficient large, such that $m^{-1} < x < m + m^{-1}$, we have

$$\frac{d}{dv} K_m(x; v) \begin{cases} > 0 & m^{-1} \leq v < x - m^{-1}, \\ < 0 & m \geq v > x - m^{-1} + \rho, \end{cases} \tag{25}$$

where $\rho$ is some small positive number. Also, note that $\frac{d^2}{dv^2} K_m(x; v) < 0$ for all $x > 0$ and $m^{-1} \leq v \leq m$. Thus, the first order derivative changes from positive to negative as $v$ changes from $m^{-1}$ to $m$ for given $x$ and sufficient large $m$, and hence, there exists $m_0$ such that $K(x; v, m^{-1})$ is increasing as a function of $v$ when $v \leq m_0$ and decreasing when $v > m_0$. For sufficient large $m$,

$$e^{-xm} \left[ \sum_{t=0}^{[m^2]} \frac{(xm)^t}{t!} - 1 - \frac{(xm)^{[m_0]+1}}{([m_0]+1)!} \right]$$

$$\leq e^{-xm} \left[ \int_{1}^{m^2} \frac{(xm)^{vm}}{\Gamma(vm+1)} d(vm) \right]$$

$$\leq e^{-xm} \left[ \sum_{t=0}^{[m^2]} \frac{(xm)^t}{t!} - 1 + \frac{(xm)^{m_0}}{([m_0]-1)!} \right], \tag{26}$$

where $[z]$ stands for the largest integer less than or equal to $z$. Using the expression for the remainder of Taylor's series, we have the LHS of (26) at least

$$1 - \frac{\frac{(xm)^{[m^2]+1}}{([m^2]+1)!} e^{x^*m}}{e^{xm}} - \frac{1}{e^{xm}} - \frac{\frac{(xm)^{[m_0]+1}}{([m_0]+1)!}}{e^{xm}}, \tag{27}$$

where $x^* \in (0, x)$. It is obvious that the expression in (27) tends to 1 as $m \to \infty$. Similarly, we have that the RHS of (26) tends to 1 as $m \to \infty$. Hence,

$$\int_{m^{-1}}^{m} K(x; v, u) dv = e^{-xm} \int_{1}^{m^2} \frac{(xm)^{vm}}{\Gamma(vm+1)} d(vm) \to 1 \quad \text{as } m \to \infty,$$

that is, Condition C1 is satisfied.

From above, we also know that Condition C3 is satisfied, since $K_m(x; v) > 0$ for all $v \in A_m$ and $x \in \mathfrak{X}$.

To verify Condition C2, for any $\delta > 0$ and $x \in \mathfrak{X}$, we want

$$\int_{|x-v|>\delta, v \in A_m} \left| K_m(x, v) \right| dv = \int_{|x-v|>\delta, v \in A_m} \frac{e^{-xm}(xm)^{vm}}{\Gamma(vm+1)} dv \to 0,$$

as $m \to \infty$. We show that for any $\delta > 0$,

$$m \sup_{|x-v|>\delta, v \in A_m} \frac{e^{-xm}(xm)^{vm}}{\Gamma(vm+1)} \to 0 \quad \text{as } m \to \infty,$$

which is equivalent to showing that

$$\log m + \log \frac{e^{-xm}(xm)^{vm}}{\Gamma(vm+1)} \to -\infty \quad \text{for all } v \in A_m, |x-v| > \delta.$$

For any $v$ such that $v \in A_m, |x - v| > \delta$, we have by Stirling's inequality for factorials,

$$\log m + \log \frac{e^{-xm}(xm)^{vm}}{\Gamma(vm+1)}$$

$$\leq \log m + \log \frac{e^{-xm}(xm)^{vm}}{[vm]!}$$

$$\leq \log m + vm \log(xm) - xm - vm \log vm + vm$$

$$= \log m + \{1 + \log(x/v) - x/v\}vm \to -\infty,$$

as $m \to \infty$, since for any given $x$ and $\delta$, there exists $q < 0$ such that $1 + \log(x/v) - x/v < q$ for all the $v \in A_m, |x - v| > \delta$.

Thus Conditions C1–C3 in Lemma 7 are all satisfied and we have that $f_m(x) \to f_0(x)$ as $m \to \infty$ for each $x > 0$. $\square$

**Lemma 8.** *Let $K_m(x; \alpha)$ be defined as in Section 12. If Condition B7\* is satisfied, then there exists a function $0 < C(x) < 1$ such that*

$$C(x) \leq \begin{cases} \int_{m^{-1} \vee x}^{x+\delta} K_m(x; mv+1)dv, & m^{-1} < x < 1, \\ \int_{x-\delta}^{m \wedge x} K_m(x; mv+1)dv, & 1 \leq x \leq m + m^{-1}, \end{cases} \tag{28}$$

*and $\int \log \frac{1}{C(x)} f_0(x)dx < \infty$.*

PROOF. For $m^{-1} < x < 1$, applying Stirling's inequality and noting that $v < x + \delta < 1 + \delta$ in the following integral, it follows that

$$\int_{m^{-1} \vee x}^{x+\delta} K_m(x; mv+1)dv$$

$$= \int_{m^{-1} \vee x}^{x+\delta} \frac{m^{mv+1}x^{mv}e^{-mx}}{\Gamma(mv+1)}dv$$

$$\geq \int_{m^{-1} \vee x}^{x+\delta} \frac{m^{mv+1}x^{mv}e^{-mx}}{\sqrt{2\pi}(mv+1)^{mv+1/2}\exp\{-(mv+1)+(12x)^{-1}\}}dv$$

$$= \sqrt{\frac{m}{2\pi}}\exp(1-(12x)^{-1})\int_{m^{-1} \vee x}^{x+\delta} \frac{x^{mv}e^{m(v-x)}}{(v+m^{-1})^{mv+1/2}}dv$$

$$\geq \frac{\sqrt{m}}{\sqrt{2\pi(1+\delta+m^{-1})}}\exp(1-(12x)^{-1})\int_{m^{-1} \vee x}^{x+\delta} \frac{x^{mv}e^{m(v-x)}}{(v+m^{-1})^{mv}}dv. \tag{29}$$

Note that

$$\int_{m^{-1}\vee x}^{x+\delta} \frac{x^{mv}e^{m(v-x)}}{(v+m^{-1})^{mv}}dv$$

$$= \int_{m^{-1}\vee x}^{x+\delta} \exp\left[mv\left\{\log\frac{x}{v+m^{-1}} - \left(\frac{x}{v}-1\right)\right\}\right]dv$$

$$= \int_{m^{-1}\vee x}^{x+\delta} \exp\left[mv\left\{\log\frac{x}{v+m^{-1}} - \left(\frac{x}{v+m^{-1}}-1\right)\right.\right.$$
$$\left.\left. +\left(\frac{x}{v+m^{-1}}-1\right) - \left(\frac{x}{v}-1\right)\right\}\right]dv$$

$$> \int_{m^{-1}\vee x}^{x+\delta} \exp\left[mv\left\{-\frac{1}{2\frac{x}{v+m^{-1}}}\left(\frac{x}{v+m^{-1}}-1\right)^2 + \frac{-x/m}{v(v+m^{-1})}\right\}\right]dv$$

$$= \int_{m^{-1}\vee x}^{x+\delta} \exp\left[\frac{-mv(x-v-m^{-1})^2 - 2x^2}{2x(v+m^{-1})}\right]dv.$$

The above inequality holds, because of that, for $0 < u < 1$,

$$\log u - (u-1) = -(1-u)^2\left\{\frac{1}{2} + \frac{(1-u)}{3} + \frac{(1-u)^2}{4} + \cdots\right\}$$

$$\geq -\frac{(1-u)^2}{2}\left\{1 + (1-u) + (1-u)^2 + \cdots\right\} = -\frac{(1-u)^2}{2u}.$$

Since $1 + \delta > x + \delta > v > x$ in the following integral, we have that

$$\int_{m^{-1}\vee x}^{x+\delta} \exp\left(\frac{-mv(x-v-m^{-1})^2 - 2x^2}{2x(v+m^{-1})}\right)dv$$

$$\geq \int_{(m^{-1}\vee x)+m^{-1}}^{x+\delta+m^{-1}} \exp\left(\frac{-m(1+\delta)(x-\tilde{v})^2 - 2x^2}{2x^2}\right)d\tilde{v}$$

$$= \sqrt{\frac{2\pi}{m}}\frac{x}{\sqrt{1+\delta}}e^{-1}\left\{\Phi\left(\frac{\delta+m^{-1}}{x/\sqrt{m(1+\delta)}}\right) - \Phi\left(\frac{m^{-1}}{x/\sqrt{m(1+\delta)}}\right)\right\}$$

$$\geq \sqrt{\frac{2\pi}{m}}\frac{x}{\sqrt{1+\delta}}e^{-1}\left\{\Phi\left(\frac{\delta+m^{-1}}{x/\sqrt{m(1+\delta)}}\right) - \Phi\left(\frac{m^{-1}}{x/\sqrt{m(1+\delta)}}\right)\right\}, (30)$$

where $\tilde{v} = v + m^{-1}$ and $\Phi(\cdot)$ is the c.d.f. of the standard normal distribution. For $m$ large, such that $\delta > m^{-1/2}$,

$$\Phi\left(\frac{\delta+m^{-1}}{x/\sqrt{m(1+\delta)}}\right) - \Phi\left(\frac{m^{-1}}{x/\sqrt{m(1+\delta)}}\right)$$

$$= \Phi\left(\sqrt{1+\delta}\,\frac{m^{1/2}\delta+m^{-1/2}}{x}\right) - \Phi\left(\sqrt{1+\delta}\,\frac{m^{-1/2}}{x}\right)$$

$$\geq \Phi(2\sqrt{1+\delta}\,\sqrt{\delta}/x) - \Phi(\sqrt{1+\delta}\,\delta/x)$$

$$\geq \Phi(2\sqrt{1+\delta}\,\delta/x) - \Phi(\sqrt{1+\delta}\,\delta/x). \tag{31}$$

The last inequality holds since we chose $\delta < 1$. Now for $u > 0$,

$$\frac{\frac{1+u^2}{u}\phi(u)}{\frac{1}{2u}\phi(2u)} = 2(1+u^2)e^{3u^2/2} \geq 2,$$

where $\phi(x) = (2\pi)^{-1/2}e^{-x^2/2}$ is the standard normal p.d.f. By the fact that

$$\frac{x}{1+x^2}\phi(x) < 1 - \Phi(x) < \frac{\phi(x)}{x}, \tag{32}$$

we have that

$$\Phi(2u) - \Phi(u) \geq \frac{1+u^2}{u}\phi(u) - \frac{1}{2u}\phi(2u) \geq \frac{1}{2u}\phi(2u).$$

Hence, the RHS of (31) is greater than

$$\frac{x}{2\delta\sqrt{2\pi(1+\delta)}}\exp\left(-\frac{2(1+\delta)\delta^2}{x^2}\right). \tag{33}$$

Now, combining the expressions (29), (30) and (33), it follows that

$$C(x) = \frac{x^2}{2\delta(1+\delta)\sqrt{2+\delta}}\exp\left(-\frac{1}{12x} - \frac{2(1+\delta)\delta^2}{x^2}\right), \qquad 0 < x < 1, \tag{34}$$

satisfies (28) for $m^{-1} < x < 1$.

Now let $m + m^{-1} > x \geq 1$. Applying Stirling's inequality, we have that

$$\int_{(x-\delta)}^{m\wedge x} K_m(x; mv+1)dv$$

$$= \int_{(x-\delta)}^{m\wedge x} \frac{m^{mv+1}x^{mv}e^{-mx}}{\Gamma(mv+1)}dv$$

$$\geq \int_{(x-\delta)}^{m\wedge x} \frac{m^{mv+1}x^{mv}e^{-mx}}{\sqrt{2\pi}(mv+1)^{mv+1/2}\exp[-(mv+1)+(12x)^{-1}]}dv$$

$$= \sqrt{\frac{m}{2\pi}}\, e^{1-(12x)^{-1}} \int_{(x-\delta)}^{m\wedge x} \frac{x^{mv}e^{m(v-x)}}{(v+m^{-1})^{mv+1/2}}dv$$

$$\geq \frac{\sqrt{m}\, e^{1-(12x)^{-1}}}{\sqrt{2\pi(x+\delta)}} \int_{x-\delta}^{x\wedge m} \exp\left[mv\left\{\log\left(\frac{x}{v+m^{-1}}\right) + \left(1-\frac{x}{v}\right)\right\}\right]dv, \tag{35}$$

since $v + m^{-1} < x + \delta$, when $m > \delta^{-1}$. Note that

$$\log u - (u-1) = (u-1)^2\left\{-\frac{1}{2} + \frac{(u-1)}{3} - \frac{(u-1)^2}{4} + \cdots\right\}$$

$$\geq (u-1)^2\left\{-\frac{1}{2} + (u-1) - (u-1)^2 + \cdots\right\}$$

$$= (u-1)^2\left\{-\frac{1}{2} - \left[(1-u) + (1-u)^2 + (1-u)^3 + \cdots\right]\right\}$$

$$= (u-1)^2\left(-\frac{1}{2} - \frac{1-u}{u}\right),$$

for $0 < u < 1$. Further, $\log u - (u-1) \geq -(u-1)^2/2$, for $1 \leq u < 2$, since $\frac{(u-1)}{3} - \frac{(u-1)^2}{4} + \frac{(u-1)^3}{5} - \cdots \geq 0$. Note that $0 < \frac{x}{v+m^{-1}} \leq \frac{1}{1-\delta}$, where $\delta < \frac{1}{2}$ without loss of generality. Now it follows that

$$\log\left(\frac{x}{v+m^{-1}}\right) + 1 - \frac{x}{v}$$
$$= \log\left(\frac{x}{v+m^{-1}}\right) - \left(\frac{x}{v+m^{-1}} - 1\right) + \left(\frac{x}{v+m^{-1}} - 1\right) - \left(\frac{x}{v} - 1\right)$$
$$\geq \left(-\frac{1}{2} - \frac{v+m^{-1}-x}{x}\right)\left(\frac{x}{v+m^{-1}} - 1\right)^2 + \frac{x}{(v+m^{-1})mv}. \quad (36)$$

Letting $\tilde{v}$ denote $v+m^{-1}$, the RHS of (36) is equal to

$$\frac{(x-2\tilde{v})(x-\tilde{v})^2}{2x\tilde{v}} + \frac{x}{mv\tilde{v}} \geq -\frac{(x-\tilde{v})^2}{2\tilde{v}^2} + \frac{x}{mv\tilde{v}},$$

since $\frac{x-2\tilde{v}}{x} \geq -1$ for $\tilde{v} < x+m^{-1}$ (i.e. $v < x$) and $x > 1$. Now,

$$\int_{x-\delta}^{x \wedge m} \exp\left[mv\left\{\log\left(\frac{x}{v+m^{-1}}\right) + \left(1 - \frac{x}{v}\right)\right\}\right] dv$$
$$= \int_{x-\delta+m^{-1}}^{(x \wedge m)+m^{-1}} \exp\left[-\frac{mv(x-\tilde{v})^2}{2(x-\delta)^2} + \frac{x}{\tilde{v}}\right] d\tilde{v}$$
$$\geq e^{1/2} \int_{x-\delta+m^{-1}}^{(x \wedge m)+m^{-1}} \exp\left[-\frac{mx(x-\tilde{v})^2}{2(x-\delta)^2}\right] d\tilde{v}$$
$$\geq e^{1/2}\sqrt{\frac{2\pi}{mx}}(x-\delta)\left\{\Phi\left(\frac{\sqrt{mx}\,(\delta - m^{-1})}{x-\delta}\right) - \frac{1}{2}\right\}$$
$$\geq e^{1/2}\sqrt{\frac{2\pi}{mx}}(x-\delta)\left\{\Phi\left(\frac{\sqrt{x}\,\delta}{2(x-\delta)}\right) - \frac{1}{2}\right\}$$
$$\geq \frac{1}{2}\delta\sqrt{\frac{e}{m}}\exp\left(-\frac{x\delta^2}{8(x-\delta)^2}\right), \quad (37)$$

for $m/2 > \delta^{-1}$, since $\Phi(z) - 1/2 > \phi(z)z$ for any $z > 0$ and since $x > v > 1 - \delta$ for all $x > 1$. Combining expressions (35) and (37) and simplifying, we conclude that

$$C(x) = \frac{\delta \exp(3/2 - 12x^{-1})}{2\sqrt{2\pi(x+\delta)}}\exp\left(-\frac{x\delta^2}{8(x-\delta)^2}\right), \qquad x \geq 1, \quad (38)$$

satisfies (28) for $1 \leq x < m + m^{-1}$.

Now for $C(x)$ defined by (34) and (38) satisfies (28). Further, by straightforward calculations, $\int \log\frac{1}{C(x)}f_0(x)dx < \infty$ under condition B7*. $\qquad\square$